\documentclass{amsart}
\usepackage{amsmath,amsthm,amssymb,amsfonts}
\usepackage[colorlinks=true]{hyperref}

\hypersetup{urlcolor=blue, citecolor=blue}

\def\arXiv#1{   {\href{http://arxiv.org/abs/#1}
   {{\mdseries\ttfamily arXiv:#1}}}}

\def\doi#1{   {\href{http://dx.doi.org/#1}
   {{\mdseries\ttfamily DOI}}}}

\newcommand{\al}{\alpha}    \newcommand{\be}{\beta}
\newcommand{\de}{\delta}    
  \newcommand{\ep}{\varepsilon}

\newcommand{\om}{\omega}    
    
\newcommand{\R}{\mathbb{R}}\newcommand{\Z}{\mathbb{Z}}

\newcommand{\pt}{\partial_t}\newcommand{\pa}{\partial}
\newcommand{\D}{{\mathrm D}}\newcommand{\les}{{\lesssim}}
\newcommand{\beeq}{\begin{equation}}\newcommand{\eneq}{\end{equation}}

\newcommand{\Sp}{{\mathbb S}}\def\CO{\mathcal {O}}

\newcommand{\supp}{\text{supp}}

\newenvironment{prf}{\noindent {\bf Proof.} }{\endprf\par}
\def \endprf{\hfill  {\vrule height6pt width6pt depth0pt}\medskip}

\numberwithin{equation}{section}

\newcommand{\gm}{\mathfrak{g}}

\def\O{{\mathcal{O}}}

\def\<{\langle}             \def\>{\rangle}
\def\({\left(}                 \def\){\right)}

\newtheorem{thm}{Theorem}[section]
\newtheorem{prop}[thm]{Proposition}
\newtheorem{coro}[thm]{Corollary}
\newtheorem{lem}[thm]{Lemma}

\theoremstyle{remark}
\newtheorem{rem}{Remark}

\theoremstyle{definition}

\title[Long time existence for SLW on space-times]
{Long time existence for semilinear wave equations on asymptotically flat space-times}

\author{Chengbo Wang}
\address{School of Mathematical Sciences\\                Zhejiang University\\                Hangzhou 310027, China}
\email{wangcbo@gmail.com}
\urladdr{http://www.math.zju.edu.cn/wang}
\date{} 

\thanks{
The author was supported in part by NSFC 11301478 and National Support Program for Young Top-Notch Talents.
}
\dedicatory{} \commby{}

\begin{document}

\begin{abstract}
We study the long time existence of solutions to nonlinear wave equations with power-type nonlinearity (of order $p$) and small data, on a large class of $(1+n)$-dimensional nonstationary asymptotically flat backgrounds, which include the Schwarzschild and Kerr black hole space-times. Under the assumption that uniform energy bounds and a weak form of local energy estimates hold forward in time, we give lower bounds of the lifespan when $n=3, 4$ and $p$ is not bigger than the critical one.
The lower bounds for three dimensional subcritical and four dimensional critical cases are sharp in general. 
For the most delicate three dimensional critical case,
we obtain the first existence result up to $\exp(c\ep^{-2\sqrt{2}})$, for
many space-times including the  nontrapping exterior domain, nontrapping asymptotically Euclidean space and Schwarzschild space-time. 
\end{abstract}



\maketitle 

\section{Introduction}
\label{sec-Intro}

Recently, there have been many advances in the understanding of the analog of the Strauss conjecture on various $(1+n)$-dimensional space-times, including exterior domains, asymptotically Euclidean manifolds, and Schwarzschild/Kerr black hole space-times. 
The problem, initiated in \cite{John79} for the Minkowski space-time, is about the small data problem for the nonlinear wave equations $\Box_g u=|u|^p$. In particular, when $p$ is greater than the Strauss exponent $p_c(n)$ \cite{Strauss81}, which is the positive root of $(n-1)p^2-(n+1)p-2=0$ ($p_c(3)=1+\sqrt{2}$, $p_c(4)=2$), the low dimensional global existence has been well understood in many settings. See \cite{DMSZ, HMSSZ, SmSoWa12} for exterior domains ($2\le n\le 4$), \cite{SoWa10, WaYu11} for asymptotically Euclidean manifolds ($n=3,4$), \cite{DaRo05, BlSt06, LMSTW} for Schwarzschild/Kerr black hole space-times and the recent work \cite{MW} for general asymptotically flat space-times.

In contrast, less is known (nonexistence of global solutions, estimates of the sharp lifespan from above and below) for the problem with $1<p\le p_c(n)$. 
 In this paper, we are interested in the investigation of the existence portion of the problem, that is, to obtain lower bounds for the lifespan, on a large class of asymptotically flat space-times, including Schwarzschild/Kerr black hole space-times. For more details on the history of the Strauss  conjecture,
see \cite{WaYu11p}. 

This problem has been relatively well understood in the setting of the exterior domains with compact  metric perturbation. Concerning the nonexistence of global solutions and the upper bound of the lifespan, see \cite{ZhouHan1} for $1<p<p_c(n)$ and 
\cite{LaZh14p} for $p=p_c(3)$. 
On the other hand, the lower bounds are known in \cite{DuZh08, Yu11, ZhZh15} for star-shaped obstacles ($2\le p\le p_c(n)$, $n=3,4$),
which agree with the upper bound, for $2\le p<p_c(3)$.

For the corresponding problem on asymptotically Euclidean manifolds, it is known to admit a weak lower bound \cite{WaYu11} for $2\le p<p_c(3)$, compared to the free case. For the Schwarzschild black hole background, the problem with $1<p<p_c(3)$ does not admit global existence in a weak sense  \cite{CG06}. 

Let $n\ge 3$, $2\le p\le p_c(n)$, for $(1+n)$-dimensional smooth asymptotically flat Lorentzian space-times $(M,g)$, we are interested in studying the long time existence for the  following nonlinear wave equations with small initial data
\begin{equation}
\label{eq-nlw}
\left\{\begin{array}{l}\Box_g u =F_p(u),\ (t,x)\in M\\
u(0,x)=u_0(x), \pt u(0,x)=u_1(x)\ .
\end{array}\right.
\end{equation}
where $F_p$ behaves like $|u|^p$ for small $u$, that is, for $|u|\ll 1$,
$$|\pa^k_u F_p(u)|\les |u|^{p-k}, k\le 2\ .$$
Typical examples include $\pm |u|^{p-1}u, \pm |u|^p$.

Comparing to the global existence of the problem with $p>p_c(n)$, the problem of existence with sharp lower bound for the subcritical case $p< p_c(n)$ is more difficult to handle, in the sense that we can not admit any loss for the weights in the key estimates (e.g., weighted Strichartz estimates) adjusted for the nonlinearity. For example, a global in time estimate with arbitrarily small loss for the weights has been sufficient for proving global existence in the setting of asymptotically Euclidean manifolds \cite{SoWa10, WaYu11}, but 
a local in time estimate with small loss will not be sufficient to prove long time existence with sharp lower bound of the lifespan for the subcritical case \cite{WaYu11}. The critical case $p=p_c(n)$ is even more delicate.

The first main result in this paper,
Theorems \ref{metaTheorem1},
 assert almost global existence of order $\exp(c \ep^{-2})$ for the four dimensional critical case with $F_p=u^2$, which is known to be sharp in general. For the three dimensional subcritical case, we also obtain
the long time existence, which is of the same order as the sharp lifespan in the Minkowski space-time and is sharp in general. For the  three dimensional critical case, it is more delicate and
we obtain the first existence result up to $\exp(c\ep^{-2\sqrt{2}})$, 
for non-flat manifolds with stationary and split metric, near spatial infinity.

The strategy of proof is to devise the local in time space-time estimates adapted for the problem, for the part near spatial infinity.
The (integrated) local energy estimates is then used for handling the errors coming from commutators with vector fields, and the terms in the compact region.
Such approach has been proved to be very robust in the recent investigation of linear and nonlinear wave equations on various space-times, particularly in the recent understanding of the waves on black hole space-times, see,  e.g., \cite{DaRoNotes, MTT}.
For the problem \eqref{eq-nlw} with supercritical power $p>p_c$, on various space-times, 
the global in time weighted Strichartz estimates
of \cite{FaWa, HMSSZ} have played a key role in the recent works, see \cite{HMSSZ, SoWa10, LMSTW, MW}. 

As our theorems are conditional results, it is useful to review the typical cases where such local energy assumption is satisfied.  Interesting examples including Schwarzschild  black hole backgrounds, slowly rotating Kerr  black hole backgrounds, nontrapping exterior domain, nontrapping asymptotically Euclidean manifolds, 
and asymptotically flat manifolds with small, time-dependent metric perturbation. 

\section{Notation and setup}

\subsection{Notation}
We use the convention that Greek indices $\al$, $\be$ range from $0$ to $n$, Latin indices $i$, $j$ from $1$ to $n$. For $(x^\al)=(t,x)$,  we shall use
polar coordinates $(t,x)=(t,r\omega)$ with $\omega\in\mathbb{S}^{n-1}$,
and the full space-time gradient $\partial=(\pt,\nabla_x)=(\partial_0,\partial_1,\cdots ,\partial_n)$. We will also use $a$, $b$ to denote multi-indices and follow the convention of summing over repeated upper and lower indices. 

With rotational vector fields $\Omega_{ij}=x^i \pa_j-x^j\pa_i$, $1\le i<j\le n$, the vector fields will be labeled as $Y=(\nabla_x,\Omega)$, $Z=(\pa, \Omega)=(\pt, Y)$.
We will also use $A\les B$ to stand for $A\le C B$ where the constant $C$ may change from line to line.
Also, $\D=\sqrt{-\Delta}$.

Let $M$ be either $M_0=\R^+\times \R^n$ or 
an exterior domain $M=\R^+\times (\R^n\backslash \mathcal{K})$ with $\mathcal{K}\subset B(0,R_0)$ 
for some $R_0 >0$ and has a smooth boundary.
We will use the following mixed-norm $L^{q_1}_t L^{q_2}_r L^{q_3}_\omega$,
$$\|f\|_{L^{q_1}_t L^{q_2}_r L^{q_3}_\omega(M)}=\left\|\left(\int_{R_0}^\infty \|f(t, r\omega)\|_{L^{q_3}_\omega}^{q_2} r^{n-1} dr\right)^{1/q_2}\right\|_{L^{q_1}(\{t\ge 0\})},$$ with trivial modification for the case $q_2=\infty$,
where $L^q_\omega$ is the standard Lebesgue space on the unit sphere $\Sp^{n-1}$. Occasionally, we will omit the subscripts. Also, $L^{q_1} L^{q_2}=L^{q_1}_t  L^{q_2}_r L^{q_2}_\omega$ and $L^q_T=L^q([0,T])$.

For convenience of presentation,  we use the notation \beeq\label{eq-BL} \|u\|_{\dot \ell^s_q(A)} = \|(\phi_j(x)u(t, x))\|_{ \dot \ell^s_q(A)}=\|\left(\|2^{js}\phi_j(x)u(t,x)\|_A\right)\|_{\ell^q_{j\in\Z}}, \eneq for a partition of unity subordinate to the homogeneous dyadic (spatial) annuli, $\sum_{j\in \Z}\phi^2_j(x)=1$ for $x\neq 0$. Similarly, we use $\ell^s_q(A)$ to denote the inhomogeneous version.

\subsection{Space-times}
We shall work on the smooth asymptotically flat Lorentzian space-times $(M,g)$.
The smooth Lorentzian metric, $g=  g_{\al\be} (t,x) \, d x^\al \, d x^\be$,
with inverse $g^{\al\be}$, is assumed to be  {\em asymptotically flat}. That is,
\begin{equation}\tag{H1} \label{H1}
g=m+g_1(t,r)+g_2(t,x), 
\end{equation}
where $(m_{\al\be})=Diag(-1,1,1,\cdots,1)$ is the standard Minkowski metric,  the first perturbation $g_1$ is {\em radial}, and
\beeq\tag{H1.1} \label{H1.1}
\|\partial^a g_{j,\al\be}\|_{\ell^{|a|+j-1}_1 L^\infty_{t,x}}\le C_{a,j},  j=1,2.
\eneq
By radial $g_1$, we mean that, we have
$$g_1(t,r)=\tilde g_{00}(t,r)dt^2+2\tilde g_{01}(t,r)dtdr+\tilde g_{11}(t,r)dr^2+\tilde g_{22}(t,r) r^2 d\omega^2\ .$$ 
Equipped with the metric $g$,  the d'Alembertian operator is given by
$\Box_g =\nabla^\al\pa_\al=  \sqrt{|g|}^{-1} \partial_\al g^{\al\be} \sqrt{|g|} \partial_\be$ with $|g| =-\det (g_{\al\be})$.

Notice that similar form of the metrics has been standard in many recent works, see, e.g., \cite{Ta13, MTT, Wa14, Wa15}.
Also, the metric models the black hole metrics, and includes Schwarzschild and Kerr black hole metric as particular case. The separation of the radial part from the general non-radial metric is crucial in many problems, see, e.g., \cite{Ta13, MTT, LMSTW}.

\subsection{Local energy assumption}
To obtain existence results, we need to assume further about the linear wave equations on such manifolds.
As been occurred  in many works (see, e.g., \cite{Ta13, MTT, Wa14, Wa15}), we make a
 {\em local energy assumption}:
There is $R_1> R_0$, such that for any solutions to the linear equation 
$\Box_g u=F$,
we have the uniform energy estimates and weak local energy estimates, for any $k\ge 0$,
\beeq \label{H2}\tag{H2} 
\|\pa^{\le k} u\|_{LE_g}
\les          \|u(0,\cdot)\|_{H^{k+1}}+\|\pa_t u(0,\cdot)\|_{H^{k}}+
\|\pa^{\le k} F\|_{L^1_tL^2_x}
\ ,
\eneq
with
$\|u\|_{LE_g}\equiv \| \pa u\|_{L^\infty_t L^2_x}+\|\psi_{R_1} \pa u\|_{ \ell^{-1/2}_\infty L^2_t L^2_x}+\| u\|_{\ell^{-3/2}_\infty  L^2_t L^2_x}$.
 Here and in what follows, we use $\|\pa^{\le k} u\|_A$ to stand for $\sum_{|a|\le k}\|\pa^{a} u\|_A$,
$\psi_R$ is a smooth nonnegative function 
\beeq\label{psiR}\psi_R=\left\{
\begin{array}{ll }
1 ,     & r>R+1,  \\
0,      &   r<R.
\end{array}
\right.
\eneq
Estimates of this form date back to \cite{Mo1}, and see, e.g., \cite{MTT, LMSTW} for more exhaustive history of such estimates.

Notice that when $g_{1}=g_{2}=0$, 
we do not need to introduce the cutoff function $\psi_{R}$ for the assumption to hold.
When the space-time permits trapped rays, such estimates without loss are known to be false in general. However, in many situations of interest, a related estimate with certain loss can be recovered. One possibility is to ask for a local energy estimate away from the trapping region, and this is exactly what we assume here. It is worth noting that we assume that there is no loss for the energy estimates (the first term in $LE_{g}$) and the lower order term from the local energy estimates (the third term in $LE_{g}$).
See Section \ref{sec-Discuss} for typical examples where such assumption holds.

Although stated in the current form, it turns out that certain weaker assumption suffices. For example, if we have for some $N$,
 \begin{eqnarray} & & \| \pa \pa^{\le k} u\|_{L^\infty_t L^2_x}+\|\<r\>^{-N} \pa \pa^{\le k} u\|_{ L^2_t L^2_x(r\ge R_1)}+\|\<r\>^{-N-1} \pa^{\le k} u\|_{ L^2_t L^2_x}
 \label{H2'}
 \\ & \les &  \|u(0,\cdot)\|_{H^{k+1}}+\|\pa_t u(0,\cdot)\|_{H^{k}}+\|\pa^{\le k} F\|_{L^1_tL^2_x}
 \ .\nonumber\end{eqnarray}
Then a similar argument as in the proof of Theorem \ref{thm-LE} will give us \eqref{H2}, see Proposition \ref{thm-H2'}. As we know from \cite{BoHa, SoWa10}, such inequalities are satisfied for nontrapping asymptotically Euclidean manifolds, with $R_1=0$ and $N>1/2$.

\section{Main results}
In this section, we present the main results in this paper.

Our first main result concerns the four dimensional critical case, $n=4$ and $p=p_{c}=2$,
for which we could prove almost global existence up to
$T_\ep=\exp(c\ep^{-2})$.
\begin{thm}\label{metaTheorem1}
Let $n\ge 4$, $p=2$, and assume \eqref{H1} and \eqref{H2}. Consider the problem \eqref{eq-nlw} with $F_p(u)=u^2$, then there exists $c>0$,  such that there is a solution $u$ in $M\cap ([0, T_{\ep})\times \R^4)$, for initial data which is sufficiently small (of size $\ep$), regular and has sufficient decay,
where \beeq\label{eq-life830}T_\ep=\left\{\begin{array}{ ll}
  \exp(c\ep^{-2}) ,   &   n=4\ ,\\
  \infty ,  &   n\ge 5\ .
\end{array}
\right.\eneq
\end{thm}

\begin{rem} Let $L_\ep(p)$ be the lifespan of the problem.
The four dimensional lower bound, $L_\ep(p_{c})\ge \exp(c \ep^{-2})$, is known to be sharp in general.
Actually, for Minkowski space-time $\R^{1+n}$ with $n\ge 3$ and $p=p_{c}$,
see  \cite{TaWa, ZhHan14} for the upper bound $\ln L_\ep\les  \ep^{-p(p-1)}$ with $F_p=-|u|^p$, as well as the lower bound $\ln L_\ep\gtrsim  \ep^{-p(p-1)}$ for $n\le 8$, in \cite{Zh92, LiZhou95, LdSo96}.
\end{rem}

When the spatial dimension is three, it turns out to be much more delicate.
\begin{thm}\label{metaTheorem2} Let $n=3$, and consider space-times with
\eqref{H1}-\eqref{H2} and stationary and {\em split} metric, near spatial infinity,
that is,
there exists $\tilde R_0\ge R_0$ such that  \beeq\label{eq-StaSplit}\tag{H3}
g=g_{00}(x)dt^2+g_{ij}(x)dx^i dx^j\ , r>\tilde R_0\ .\eneq 
Then the problem \eqref{eq-nlw} admit
solution up to 
 \beeq\label{eq-life1}T_\ep=\left\{\begin{array}{ ll}
  c\ep^{-\frac{p(p-1)}{1+2p-p^2}}    &    2\le p< 1+\sqrt{2}\\
  \exp(c\ep^{-2\sqrt{2}})    &   p= 1+\sqrt{2}
\end{array}
\right.\eneq
 for initial data which is sufficiently small of size $\ep$ and regular with sufficient decay. 
Moreover, for the general case when \eqref{eq-StaSplit} is not necessarily satisfied, the similar result hold with
 \beeq\label{eq-life2}T_\ep=\left\{\begin{array}{ ll}
  c\ep^{-\frac{p(p-1)}{1+2p-p^2}}    &    2\le p< 1+\sqrt{2}\\
  \exp(c\ep^{-\sqrt{2}})    &   p= 1+\sqrt{2}
\end{array}
\right.\eneq
for the problem with small, compactly supported data.
 \end{thm}

\begin{rem} Concerning the three dimensional critical case,
our lower bound $L_\ep \ge \exp(c\ep^{-2\sqrt{2}})$ 
is the best available lower bound, for any space-times, except for the well-understood Minkowski case $\ln L_\ep\simeq  \ep^{-(\sqrt{2}+2)}$ \cite{Zh92}.
\end{rem}
\begin{rem} 
For $2\le p<1+\sqrt{2}$, the lower bound of the lifespan, $L_\ep(p)\ge c \ep^{\frac{p(p-1)}{p^2-2p-1}}$, is sharp in general. In fact, when $M$ is the Minkowski space-time $\R^{1+3}$, it is known that $L_\ep(p)\simeq  \ep^{\frac{p(p-1)}{p^2-2p-1}}$, for $F_p=-|u|^p$ with $1<p<1+\sqrt{2}$, see 
\cite{Ld90, Zh92}. 
\end{rem}

For more precise statements of Theorems \ref{metaTheorem1}-\ref{metaTheorem2}, see Theorem \ref{metaTheorem1-2}, Theorem \ref{metaTheorem2.2} and Theorems \ref{metaTheorem}.

The rest of the paper is organized as follows.
In Section \ref{sec-Sobolev}, we collect weighted Sobolev type estimates, including an improved version of the trace
estimates
and a weighted variant of the trace estimates.
In Section \ref{sec-General}, for quadratic semilinear problems, we prove the high dimensional existence results, Theorem \ref{metaTheorem1}. To achieve this goal, 
we present a weighted $L^2_t L^\infty_r L^2_\omega$ estimate, and a space-time estimate adapted for the problem near spatial infinity, Theorem \ref{thm-wStri-homo}. For handling the errors coming from commutators, and the terms in the compact region, we give the desired local energy estimates for $\Box_\gm$, Theorem \ref{thm-LE}, based on the local energy estimates with small metric perturbation as well as the local energy assumption \eqref{H2}.
In Section \ref{sec-SplitMetric},
for the three dimensional space-times with stationary and split metric near spatial infinity,
we devise a space-time estimate, Theorem \ref{ap-thm-le-0wS},
and give the proof of Theorem \ref{metaTheorem2},
for the case when \eqref{eq-StaSplit} is satisfied.
The remaining part of Theorem \ref{metaTheorem2} is proved in
Section \ref{sec-compact}.
Finally, in Section \ref{sec-Discuss},  we conclude this paper with a discussion of the implications of Theorems \ref{metaTheorem1}-\ref{metaTheorem2}, for many space-times of interest. 

\section{Sobolev-type estimates}
\label{sec-Sobolev}

In this section, we give several Sobolev type estimates. At first, let us present an improved version of the trace estimates. 
\begin{lem}[Trace estimates]\label{thm-trace}
Let $n\ge 2$ and $1/2< s<n/2$, then we have
\beeq\label{eq-trace0}
\|r^{n/2-s} f\|_{L^\infty_r H^{s-1/2}_\omega}\le C_{s} \|f\|_{\dot H^s}, \ 
\|r^{(n-1)/2} f\|_{L^\infty_r L^{2}_\omega}\les \|f\|_{\dot B^{1/2}_{2,1}},
\eneq
\beeq\label{eq-trace}
\|r^{n/2-s} f\|_{L^\infty_r H^{s-1/2-}_\omega}\le C_{s,\de} \|f\|_{\dot B^s_{2,\infty}}, \ 
 \|f\|_{\dot B^{-s}_{2,1}}
 \les \|r^{s-n/2} f\|_{L^1_r H^{-s+1/2+}_\omega}\ , 
\eneq
\beeq\label{eq-trace2-sd}
\|r^{(n-1)(1/2-1/p)} f\|_{L^p_r L^{2}_\omega}\les \|f\|_{\dot H^{1/2-1/p}}, 2\le p<\infty\ ,
\eneq for any $f\in C_0^\infty (\R^n)$.
Here $\dot B^{s}_{p,q}$ is the homogeneous Besov space, $H^s_\omega$ is the Sobolev space on the unit sphere, and  $a\pm$ stands for $a\pm \de$ for any sufficiently small $\de> 0$.
\end{lem}
\begin{prf}
The estimate \eqref{eq-trace0} has been well-known, see, e.g., (1.3), (1.7) in \cite{FaWa} and references therein.
The inequality \eqref{eq-trace2-sd} is due to \cite{LiZhou95}, see also \cite{HWY3} for alternative proof using real interpolation and \eqref{eq-trace0}.
 To get \eqref{eq-trace} with $\de>0$, for given $f$, it suffices to apply the known estimates for Littlewood-Paley projection $S_j f$ with $s\pm \de$ and sum over to prove the estimates for small enough $\de>0$, here $S_j=\phi_j*$ denotes the standard Littlewood-Paley projection operator, with
$\phi_j(x)=2^{jn}\phi(2^j x)$,
 $\supp\ \hat \phi\subset \{|\xi|\in [2^{-2}, 2^{2}]\}$.
  More precisely, for given $s\in (1/2, n/2)$, let $\de>0$ such that $s\pm \de\in  (1/2, n/2)$. Then for the part with $j\le 0$,
$$
\|\sum_{j\le 0} S_j f (\omega)\|_{ H^{s+\de-1/2}_\om}
\le\sum_{j\le 0}\| S_j f (\omega)\|_{ H^{s+\de-1/2}_\om}
\le C \sum_{j\le 0}\| S_j f \|_{\dot H^{s+\de}}
\le C \| f \|_{\dot B^{s}_{2,\infty}}\ .
$$
The counterpart with $j>0$ follows similarly, which gives \eqref{eq-trace} for $r=1$ and then for any $r>0$ by rescaling. 
\end{prf}

We shall also exploit the following weighted variant of the trace estimates. 
\begin{lem}[Weighted trace estimates]\label{thm-trace-w}
Let $n\ge 3$, $\al\in (1/2, 1]$ and $\be \in (\al-n/2, n/2)$. Then we have
\beeq\label{eq-trace-w}
\|r^{n/2-\al+\be}u\|_{L^\infty_{r} H^{\al-1/2-}_\omega}\les 
\|r^{\be} \D^{\al} u\|_{L^{2}}\ .
\eneq
\end{lem}
\begin{prf}
At first, we observe that the proof
of \eqref{eq-trace-w1} can be reduced to
 the following weak estimates
\beeq\label{eq-trace-w1}
\|u(\om)\|_{H^{\al-1/2}_\omega(\Sp^{n-1})}\les 
\|r^{\be} \nabla u\|_{L^{2}}^{\al}\|r^{\be}  u\|_{L^{2}}^{1-\al}\ , \al\in [1/2,1],\ \be \in (\al-n/2, n/2)\ .
\eneq
In fact, the estimate \eqref{eq-trace-w} with $\al=1$ follows directly from \eqref{eq-trace-w1} by rescaling.
When $\al\in (1/2, 1)$,
applying the  estimate \eqref{eq-trace-w1} for  Littlewood-Paley projection $S_j u$ and summation gives us
\beeq
\|u(\om)\|_{H^{\al-1/2}_\omega(\Sp^{n-1})}\les 
\|r^{\be} 2^{j\al} S_{j} u\|_{\ell_{j}^{1}L^{2}}\equiv 
\|u\|_{\dot B^{\al,\be}_{2,1}}
\ , \al\in [1/2,1]\ .
\eneq
Applying similar proof as in Lemma \ref{thm-trace} and scaling, we get
\beeq\label{eq-trace-w3}
\|r^{n/2-\al+\be}u\|_{L^\infty_{r} H^{\al-1/2-}_\omega}\les 
\|r^{\be} 2^{j\al} S_{j} u\|_{\ell_{j}^{2}L^{2}}
\ , \al\in (1/2,1)\ ,
\eneq
which gives us the desired
 \eqref{eq-trace-w} with $\al\in (1/2,1)$, provided that we have
\beeq\label{eq-trace-w4}
\|r^{\be} 2^{j\al} S_{j} u\|_{\ell_{j}^{2}L^{2}}
\simeq \|r^{\be} \D^{\al} u\|_{L^{2}}
\ , \al\in (1/2,1)\ ,  \be \in (\al-n/2, n/2)\ .
\eneq

It remains to prove \eqref{eq-trace-w1} and \eqref{eq-trace-w4}.
To prove \eqref{eq-trace-w4}, we recall some well known facts from harmonic analysis.
It is known (see, e.g., \cite[Theorem 7.21]{Sch13})
 that if $T$ is a strong Calder\'on-Zygmund operator, then
\beeq\|T(f)\|_{L^p(wdx)}\le C\|f\|_{L^p(wdx)}, w\in A_p, p\in (1,\infty)\ ,\eneq
where $A_p$ denotes the Muckenhoupt $A_p$ class of weight functions.
Based on this fact, it is easy to adapt the argument in, e.g., \cite[Section 8.2]{Sch13} to conclude the following weighted Littlewood-Paley square-function estimate 
\beeq\label{eq-LPw}\|w S_j f\|_{L^p \ell^2_j}\simeq \|w f\|_{L^p}, w^{p}\in A_p, f\in L^p(wdx), p\in (1,\infty)\ .\eneq
As
$r^{2\be}\in A_{2}$ if and only if $|\be|<n/2$ (see, e.g., \cite[Example 7.1.7]{Gra14}),
we get \eqref{eq-trace-w4} from \eqref{eq-LPw}.

For the proof of \eqref{eq-trace-w1}, we
recall the weighted Hardy-Littlewood-Sobolev estimates of Stein-Weiss 
$$\|r^{\be-\al} u\|_{2}\les \|r^{\be}\D^{\al} u\|_{2}, \al\in (0,n/2), \be\in (\al-n/2,n/2)\ .$$ Then for any $n\ge 3$, $\al\in (0,1]$ and $\be\in (\al-n/2,n/2)$, we have
$$\|r^{\be-\al} u\|_{2}\les \|r^{\be}\D^{\al} u\|_{2}
\simeq
\|r^{\be} 2^{j\al}S_{j} u\|_{\ell^{2}_{j}L^2}
\les \|r^{\be}2^{j}S_{j} u\|_{2}^{\al}\|r^{\be} u\|_{2}^{1-\al}
\les \|r^{\be}\nabla u\|_{2}^{\al}\|r^{\be} u\|_{2}^{1-\al}
\ .$$

Let $\phi$ be a cutoff function of $B_{2}\backslash B_{1/2}$ which equals one for $|x|=1$, we get from
\eqref{eq-trace0} that for $\al \in [1/2, 1]$,
\begin{eqnarray*}\|u(\omega)\|_{H^{\al-1/2}_\omega}&\les &
\| \nabla (\phi u)\|_{L^{2}}^{\al}\|\phi u\|_{L^{2}}^{1-\al}\\
&\les& \|r^{\be} \nabla u\|_{L^{2}}^{\al}\|r^{\be}  u\|_{L^{2}}^{1-\al}
+\|r^{\be-\al} u\|_{L^{2}}\\
&\les& \|r^{\be} \nabla u\|_{L^{2}}^{\al}\|r^{\be}  u\|_{L^{2}}^{1-\al}\ ,
\end{eqnarray*}
which gives us \eqref{eq-trace-w1}.
\end{prf}

We will also need the following variant of the Sobolev embeddings.
\begin{lem}[Weighted Sobolev estimates]\label{thm-Sobo}
Let $n\ge 2$, $R\ge 3$ and $k=[(n+2)/2]$ be the integral part of $(n+2)/2$. We have
\beeq\label{eq-Sobo}
\|r^{b}u\|_{L^q_{r} L^\infty_\omega(r\ge R+1)}\les  \|r^{b-(n-1)/p+(n-1)/q} Y^{\le k} u\|_{L^p_{r}L^{2}_\omega(r\ge R)}
\eneq
for any $b\in\R$, $2\le p\le q\le \infty$. Moreover, for any $2\le p\le q\le {q_2}<\infty$, $b\in\R$, $s\in \Z$ with
$s\ge n/2-n/{q_2}$, we have
\beeq\label{eq-Sobo2}
\|r^{b}u\|_{L^q_{r} L^{q_2}_\omega(r\ge R+1)}\les  \|r^{b-(n-1)/p+(n-1)/q} Y^{\le s} u\|_{L^p_{r}L^{2}_\omega(r\ge R)}\ .
\eneq
\end{lem}
It is essentially proved in  \cite[Lemma 3.1]{LMSTW}, in the case of $n=3$ and particular choices of $p$ and $q$, based on Sobolev's lemma in polar coordinates $(r,\omega)$.

\section{High dimensional existence}\label{sec-General}
In this section, we give the proof of Theorem \ref{metaTheorem1}, for which we record a more precise version here:
\begin{thm}\label{metaTheorem1-2}
Let $n\ge 4$, 
$k=[(n+2)/2]$ be the integral part of $(n+2)/2$. 
 Considering the problem \eqref{eq-nlw} with $p=2$ and $F_p=u^2$, posed on space-time $(M,g)$ with \eqref{H1}-\eqref{H2}, then there exists $R, c>0$ and $\ep_0>0$, such that 
for $T\le T_\ep$ and any $(u_0, u_1)$ with \beeq\label{9-Th-as}\|Y^{\le k} u_0\|_{H^1}+\|Y^{\le k}u_1\|_{L^2}+
\|\psi_R Y^{\le k} u_1\|_{\dot H^{-1}}=\ep\le \ep_0\ ,\eneq
 there is an unique solution $u$ in $M\cap ([0, T_\ep)\times \R^n)$, such that
$$
 \| Z^{\le k} u\|_{LE_g}
+\|\<r\>^{-(n-1)/4} Z^{\le k} u\|_{L_T^2 L^2}+
A_T^{-1/4}\|\<r\>^{-(n-3)/4} \psi_R Z^{\le k} u\|_{L_T^{4}L^{2}}
\le C\ep \ ,$$
 where
$T_{\ep}$ is given in \eqref{eq-life830}
and $$A_T=\left\{
\begin{array}{ll}
    \ln T,   &  n= 4\ ,\\
    1  ,  &   n\ge 5 
 \ .
\end{array}
\right.$$   
\end{thm}
\begin{rem} 
In the case of $n=3$ with $p\in [2,1+\sqrt{2}]$, for any $\de\in (0,2)$, similar proof will yield results with $T_\ep=c\ep^{-\frac{p(p-1)}{1+2p-p^2+\de p}}$, provided that
$$\|Y^{\le 2} u_0\|_{H^1}+\|Y^{\le 2}u_1\|_{L^2}+
\|\psi_R Y^{\le 2} (u_0,u_1)\|_{\dot B^{1/2-(1-\de)/p}_{2,p}\times\dot B^{-1/2-(1-\de)/p}_{2,p}}=\ep\le \ep_0\ .$$
The technical reason for the introduction of $\de$ is that Lemma \ref{thm-MetTa9} is not available for $s=-1$ and $n=3$. 
\end{rem}

This result is mainly based on the following space-time estimates for the D'Alembertian with small metric perturbation.
\begin{thm}[Space-time estimates]\label{thm-wStri-homo}
Let $n\ge 4$, consider 
 \beeq\label{eq-Box-per1}
\Box_h=\pt^2-\Delta +h_0^{\al\be}(t,x)\pa_\al\pa_\be+h_1^\al(t,x)\pa_\al\eneq 
with  $h^{\al\be}_0=h_0^{\be\al}\in C^2$, $h_1\in C^{1}$, and
\beeq\label{eq-kss-as2'}\|h_0\|_{\dot \ell^0_1 L^\infty_{t,x}}+\|(\pa h_0, h_1)\|_{\dot \ell^1_1 L^\infty_{t,x}}+
\|(\pa^2 h_0, \pa h_1)\|_{\dot \ell^2_1 L^\infty_{t,x}}\le \delta_1\ .\eneq
Then there exist a $\de_1>0$ such that for any $T>2$ and $q\ge 2$, we have
\begin{eqnarray}
& & (\ln T)^{-1/q}\|\<r\>^{-1/q} w\|_{L_T^{q}L^{2}}
+\|\<r\>^{-1/q-} w\|_{L_T^{q}L^{2}}  \label{eq-9-06}  \\
& \les & 
\|\pa w(0,\cdot)\|_{\dot H^{-1}}
+\|r^{-(n-3)/2}F_{1}\|_{L^2L^1H^{-1/2+}_{\omega}}
+\|F_{2}\|_{\dot \ell^{1/2}_1 L^2\dot H^{-1}+\dot \ell^{3/2}_2 L^2L^2+L^1\dot H^{-1}}\nonumber
\end{eqnarray}
for any $w$ with $\Box_h w=F_{1}+F_{2}$.
\end{thm}
We note that the left hand side resembles the KSS type norm \cite{KSS1}, while
the right hand norm for $F_1$ resembles \cite[Theorem 6.4]{LdSo96} for $n\ge 3$ and \cite[Lemma 3.5]{ZhZh15} for $n=4$.

\subsection{Proof of Theorem \ref{thm-wStri-homo}}\label{sec:5.2}
To prove the space-time estimates \eqref{eq-9-06}, we need to exploit the following microlocal version of localized energy estimates. 
\begin{lem}[\cite{MetTa07}, Theorem 1]\label{thm-MetTa9}
Let $n\ge 2$, $|s|<(n-1)/2$ and $|s|\le 1$. Under the same assumptions as in  Theorem \ref{thm-wStri-homo},
we have
\begin{equation}
\label{eq-MetTa9}
\|\pa w\|_{L^\infty_t \dot H^s\cap X^s}\les \|\pa w(0)\|_{\dot H^s}+\|\Box_h w\|_{L^1_t \dot H^s+(X^{-s})'}\ .
\end{equation} Here $(X^{-s})'$ denotes the dual space to $X^{-s}$,
$$
\|f\|_{X^s}=\|2^{ks} S_k f\|_{l^2_k X_k}, \|f\|_{X_k}=2^{k/2}\| f\|_{L^2(A_{\le -k})}+\sup_{j> -k}\||x|^{-1/2} f\|_{L^2(A_{j})}, 
$$
with $A_j=\R_+\times\{|x|\sim 2^j\}$, $A_{\le -k}=\cup_{j\le -k} A_j$, $j, k\in \Z$.
\end{lem}

For the $X^s$ norm, it is known  that we have the following
\begin{lem}\label{thm-MetTa29}
Let $n\ge 2$ and $0<s<(n-1)/2$. Then
\begin{equation}
\label{eq-MetTa29}
\|u\|_{\dot \ell^{-1/2}_\infty L^2_{t,x} }\les \|u\|_{X^{0}
}\ ,
\|r^{-1/2-s} u\|_{L^2_{t,x} }\les \|u\|_{X^{s}
}\ .
\end{equation}
In addition, by duality, we have
\begin{equation}
\label{eq-MetTa39}
\|u\|_{(X^{s})'}\les
\|r^{1/2+s} u\|_{L^2_{t,x} }
\ .
\end{equation}
Moreover,
for any $s$,
\begin{equation}
\label{eq-MetTa49}
\|u\|_{X^{s}
}\simeq \|\nabla u\|_{X^{s-1}}
\ .
\end{equation}
\end{lem}
See \cite{Ta08, MetTa07} for \eqref{eq-MetTa29}.  

In addition, we observe that $X^{s}$ can also control certain weighted $L^2_t L^\infty_r L^2_\omega$ norm, which plays a key role in the proof.
\begin{lem}[$L^{2}L^{\infty}L^2$ estimates]\label{thm-LE-key-observ}
Let $n\ge 3$ and $s\in(1/2,(n-1)/2)$. Then for any $\de>0$, we have
\begin{equation}
\label{eq-LE-key-observ}
\|r^{(n-1)/2-s}u\|_{L^{2}_{t} L^{\infty}_{r} H^{\min(1/2, s-1/2-\de)}_{\omega} }\les \|u\|_{X^{s}
}\ .
\end{equation}
\end{lem}
\begin{prf}
Let $\al\in (1/2, 1]\cap (1/2, s)$, then 
\begin{eqnarray*}
\|r^{(n-1)/2-s}u\|_{L^{2}_{t} L^{\infty}_{r} H^{\al-1/2}_{\omega} } & \les & 
\|r^{-1/2-(s-\al)} \D^\al u\|_{L^{2}_{t,x}}
\\
 & \les & \| \D^\al u\|_{X^{s-\al}}\les \|u\|_{X^s}\ ,
\end{eqnarray*}
where we have used the weighted
trace estimates, Lemma \ref{thm-trace-w} with $\be=-1/2-(s-\al)$, in the first inequality
and \eqref{eq-MetTa29} in the second inequality.
\end{prf}

With these preparations, we could present the proof of Theorem \ref{thm-wStri-homo}.
Since $n\ge 4$, we could apply Lemmas  \ref{thm-MetTa9} and  \ref{thm-MetTa29}  with $s=-1$ to obtain
\begin{eqnarray*}
\|w\|_{L^\infty L^2 \cap \dot \ell^{-1/2}_\infty L^2_{t,x}} & \les & \|w\|_{L^\infty L^2 \cap X^0} 
\\
&\les& \|\nabla w\|_{L^\infty \dot H^{-1} \cap X^{-1}}
\\
 & \les
 & \|\pa w(0,\cdot)\|_{\dot H^{-1}}
+\|\Box_h w\|_{(X^1)'+L^1 \dot H^{-1}}\ .
\end{eqnarray*}
Recall that
Lemmas \ref{thm-MetTa29} and
\ref{thm-LE-key-observ}
 with $s=1$
 give us
\begin{eqnarray*}
&&\|u\|_{\dot \ell^{-1/2}_\infty L^2 \dot H^1 
\cap
\dot \ell^{-3/2}_2 L^2  L^2}
+\|r^{(n-3)/2}u\|_{L^2L^\infty H^{1/2-}_{\omega}}
\\&\les &
\|\nabla u\|_{\dot \ell^{-1/2}_\infty L^2 L^2}
+\|u\|_{\dot \ell^{-3/2}_2 L^2  L^2}
+\|u\|_{X^1}\\
&\les&
\|u\|_{X^{1}
}\ .\end{eqnarray*}
By duality,
we get $\|w\|_{L^\infty L^2 \cap \dot \ell^{-1/2}_\infty L^2_{t,x}}$ controlled by the right hand side \eqref{eq-9-06}.

To complete the proof, we recall 
the well-known fact that for $T\ge 2$,
\begin{eqnarray}
 &  &
T^{-\al}\|\<r\>^{-1/2+\al} u\|_{L^2_T L^2_x}+
(\ln T)^{-1/2} \|\<r\>^{-1/2} u\|_{L^2_T L^2_x}
+\|\<r\>^{-1/2-} u\|_{L^2_T L^2_x}
\label{KSS}
 \\
 & \les & 
\|u\|_{L^\infty_T L^2_x \cap  \ell^{-1/2}_\infty L^2_T L^2_x },  
\nonumber
\end{eqnarray}
which is known as KSS type estimates, and the proof is standard, see, e.g., \cite{KSS1, MetSo06, Wa14, Wa15}. As a consequence, we see that
\begin{eqnarray*}
&& (\ln T)^{-1/q}\|\<r\>^{-1/q} u\|_{L^q_T L^2_x}
+\|\<r\>^{-1/q-} u\|_{L^q_T L^2_x}  \\
& \les & (
(\ln T)^{-1/2}
\|\<r\>^{-1/2} u\|_{L^2_T L^2_x}
+\|\<r\>^{-1/2-} u\|_{L^2_T L^2_x}
)^{2/q}
\|u\|_{L^\infty_T L^2_x}^{1-2/q}\\
&\les&  \|u\|_{L^\infty_T L^2_x \cap  \ell^{-1/2}_\infty L^2_T L^2_x }\ ,
\end{eqnarray*}
 which completes the proof of \eqref{eq-9-06}.

\subsection{Proof of Theorem \ref{metaTheorem1-2}}
With help from Theorem \ref{thm-wStri-homo}, the techniques from \cite{LMSTW} essentially apply in the current situation and will 
give the proof of 
 Theorem \ref{metaTheorem1-2}.

More precisely, let
\begin{equation}
   \|\phi\|_{X_R^k} =   A_T^{-\frac 14}
  \| \psi_R \ r^{-\frac{n-3}4}Z^{\le k} \phi \|_{L^{4}_T L^{2}} +
  \| Z^{\le k} \phi\|_{LE_g} 
  +\|\<r\>^{-\frac{n-1}4} Z^{\le k} u\|_{L^2 L^2}
  ,\label{eq-9-X-norm}
\end{equation}
\beeq\label{eq-9-N-norm}
 \|F\|_{N_{R}^k} =  
 \| \psi_R^2\  r^{-\frac{n-3}2}Z^{\le k} F\|_{L_T^2L^1H^{-1/2+}_{\omega}}+\|Z^{\le k} F\|_{L_T^1L^2 L^2} \ .
\eneq
At first, the following lemma gives us the control of the solutions for the linear wave equations in $X^k_R$.
\begin{lem}\label{thm-9-keystep1}
For $n\ge 4$, there exists $R>R_1$ such that 
\begin{equation}\label{eq-9-linest}
\|u\|_{X_R^k} \lesssim \|Y^{\le k} u_0\|_{H^{1}}+\|Y^{\le k}u_1\|_{L^2}+\|\psi_R Y^{\le k}u_1\|_{\dot H^{-1}}
+\|\psi_R Z^{\le k-1}F(0,\cdot)\|_{\dot H^{-1}}
+ \|F\|_{N_{R}^k} \ ,
\end{equation} for any $k\ge 0$ and solutions to
\beeq\label{lineqn}  \Box_g u = F, \qquad  u (0,x) = u_0(x), \qquad \pt u (0,x) = u_1(x)\ .\eneq
\end{lem}
We postpone the proof to the end of this section, see Subsection \ref{sec-LE}.

Moreover, for the forcing term $F=u^2$, we have the nonlinear estimates
\begin{lem}\label{thm-9-keystep}
Let  $n\ge 4$, $k=[(n+2)/2]$, then
\beeq
\|u^2\|_{N_{R}^k}  \les   \|r^{-\frac{n-3}4} \psi_R Z^{\le k} u\|_{L^4 L^2 L^2}^2
+\|\<r\>^{-\frac{n-1}4} Z^{\le k} u\|_{L^2 L^2 L^2}^2
\les A_T^{ 1/2} \|u\|_{X_R^k}^2 \label{eq-9-keystep}\eneq
\beeq\label{eq-9-keystep2}\|u^2-v^2\|_{N_{R}^{0}}\les A_T^{1/2} (\|u\|_{X_R^k}+\|v\|_{X_R^k})\|u-v\|_{X_R^{0}}\ ,\eneq
for any interval $t\in [0,T]$ with $T\ge 2$.\end{lem}

\begin{prf}
Here, we give only the proof of \eqref{eq-9-keystep}, and omit the similar proof of \eqref{eq-9-keystep2}.
At first,  we  observe that
$$|Z^{\le k} u^2|\lesssim  \sum_{i\le j, i+j\le k}|Z^i u| |Z^{j} u|
\les
|u||Z^{\le k}u|+
\sum_{1\le i\le j, i+j\le k}|Z^i u| |Z^{j} u|\equiv I+II
\ .
$$ 

We start with bounding the first term in \eqref{eq-9-N-norm}.  By the Sobolev embeddings on $\Sp^{n-1}$ and H\"older's inequality, as $k>(n-1)/2$, we have
\beeq\label{eq-521bd1}
\| \psi_R^2\  r^{-\frac{n-3}{2}}Z^{\le k} u^2\|_{L_T^2L^1H^{-1/2+}_{\omega}}\les
\| \psi_R\  r^{-\frac{n-3}{4}}Z^{\le k} u\|_{L_T^4L^2L^{2}}^2
\les
A_T^{\frac 12} \|u\|_{X_R^k}^2\ .
\eneq

We now proceed to the second term in \eqref{eq-9-N-norm}.  We first observe
that
$$\|I \|_{L^1L^2L^2}
  \lesssim \|\<r\>^{(n-1)/4}
u\|_{L^2 L^\infty}
 \|\<r\>^{-(n-1)/4}  Z^{\le k} u\|_{L^2L^2L^{2}}
$$
Applying \eqref{eq-Sobo} and Sobolev embedding, it follows that for $k\ge [n/2]+1$,
\[
\|I \|_{L^1L^2L^2}\les
 \|\<r\>^{-(n-1)/4}  Z^{\le k} u\|_{L^2L^2L^{2}}^2
.\]
Similarly,
for any $1\le i\le j\le k-i$, there exist
$q_1, q_2\in (2,\infty)$ such that
$$k-j\ge n/{q_1}, k-i\ge n/{q_2}, 1/2=1/{q_1}+1/{q_2}$$
and then \eqref{eq-Sobo2} and Sobolev embedding give us
\begin{eqnarray*}
\||Z^i u| |Z^j u| \|_{L^1L^2L^2}&
  \lesssim &
  \|\<r\>^{(n-1)/{q_2}-(n-1)/4}
Z^i u\|_{L^2 L^{q_1}}
  \|\<r\>^{(n-1)/{q_1}-(n-1)/4}
Z^j u\|_{L^2 L^{q_2}}\\
&  \lesssim &
  \|\<r\>^{-(n-1)/4}
Z^k u\|_{L^2 L^{q_2}}^2\ .
\end{eqnarray*}
This completes the proof of Lemma \ref{thm-9-keystep}.
\end{prf}

Equipped with Lemma \ref{thm-9-keystep1}-\ref{thm-9-keystep}, it is a standard procedure to prove  
Theorem \ref{metaTheorem1}
by iteration in $X^k_R$, showing
$\|u\|_{{X}^k_R}\les \ep$ in the time interval $[0,T_\ep]$ 
with $T_{\ep}$ given in \eqref{eq-life830}.

\subsection{Proof of Lemma \ref{thm-9-keystep1}}\label{sec-LE}
In this subsection, under the assumptions \eqref{H1} and \eqref{H2},  we give the proof of 
Lemma \ref{thm-9-keystep1}
with help from Theorem \ref{thm-wStri-homo}.

At first, we record the  following classical local energy estimate, which follows by integrating $\Box u$ against a multiplier of the form $f(r)(\partial_r
 + \frac{n-1}{2r})u$ with appropriate choice of $f$,
see, e.g., \cite{MetSo06, Wa14}.
\begin{lem}\label{thm-LE0}
Let $n\ge 3$, 
 \beeq\label{eq-Box-per}
\Box_h=\pt^2-\Delta +h_0^{\al\be}(t,x)\pa_\al\pa_\be+h_1^\al(t,x)\pa_\al+h_2(t,x)\eneq with $h^{\al\be}_0=h_0^{\be\al}$, $h_0\in C^{1}(\R\times \R^n)$, $h_1, h_2\in C(\R\times \R^n)$. Then we have
\beeq\label{eq-kss-inh}\|u\|_{LE}
\les\|\pa u(0)\|_{L^2_x}+\|F\|_{L^1_t L^2_x+ \dot \ell^{1/2}_1 L^2_{t,x}}\ ,\eneq
with $\|u\|_{LE}\equiv\|\pa u\|_{L^\infty_t L^2_x}+\|\pa u\|_{\dot \ell^{-1/2}_\infty L^2_{t,x}}+\| u\|_{\dot \ell^{-3/2}_\infty L^2_{t,x}}$, provided that
\beeq\label{eq-kss-as-inh}\|h_0\|_{\dot \ell^0_1 L^\infty_{t,x}}+\|(\pa h_0, h_1)\|_{\dot \ell^1_1 L^\infty_{t,x}}+
\| h_2\|_{\dot \ell^2_1 L^\infty_{t,x}}\ll 1\ .\eneq
\end{lem}
Then,
we could strengthen the local energy assumption to estimates involving rotational vector fields. 
\begin{thm}\label{thm-LE}
Let $n\ge 3$ and assume \eqref{H1} and \eqref{H2}. Then we have
$$\|Z^{\le k} u\|_{LE_g}
\les \|Y^{\le k} u(0,\cdot)\|_{H^{1}}+\|Y^{\le k} \pt u(0,\cdot)\|_{L^{2}}
+
\|Z^{\le k} F\|_{L^1_tL^2_x} $$
for any solutions to $\Box_g u=F$
and any $k\ge 0$. 
\end{thm}
\begin{prf}
By using the equation, we need only to prove.
\beeq\|Z^{\le k} u\|_{LE_g}
\les\|Z^{\le k} \pa^{\le 1}u(0,\cdot)\|_{L^{2}}
+
\|Z^{\le k} F\|_{L^1_tL^2_x} \ .
 \label{eq-LE-2}\eneq
Let $R_2\ge R_1+1$ to be determined, and consider first $v=\psi_{R_2} u$. From the support property, by introducing $$g_0=m (1-\psi_{R_2-1})
+g \psi_{R_2-1}\ ,$$ we see that
$\Box_{g_0} v=\Box_g v$.
Thanks to \eqref{H1} with \eqref{H1.1}, we can set $R_2$ large enough such that \eqref{eq-kss-as2'}
 is satisfied for $\Box_{g_0}$, and the commutators between $Z$ and $\Box_{g_0}$ produce only small error terms comparing to the local energy norm.
That is, by Lemma \ref{thm-LE0}, we have
\beeq\label{eq-kss-inh2}
\|Z^{\le k}v\|_{LE_g}\les\|Z^{\le k}v\|_{LE}\les\|\pa Z^{\le k}v(0)\|_{L^2_x}+\|Z^{\le k}\Box_{g_0}v\|_{L^1_t L^2_x+ \dot\ell^{1/2}_1 L^2_{t,x}}\ ,\eneq
for $R_2$ large enough
(see e.g. (3.16) of \cite{Wa14} for details of similar proof).

Fixing such $R_2$, and noticing that
$\Box_{g_0} v=\Box_g v=\psi_{R_2} F+[\Box_g, \psi_{R_2}] u$,
we obtain
\begin{eqnarray*}
\|Z^{\le k} v\|_{LE_g}&\les & \|\pa Z^{\le k} v(0, \cdot)\|_{L^2}+\|Z^{\le k} F \|_{L^1_t L^2_x}+\|Z^{\le k}
[\Box_g, \psi_{R_2}] u
\|_{\dot \ell^{1/2}_1 L^2_{t,x}}\\
 & \les &  \| Z^{\le k} \pa^{\le 1} u(0, \cdot)\|_{L^2}+\|Z^{\le k} F \|_{L^1_t L^2_x}+\|\pa^{\le k+1} u
\|_{ L^2_{t,x}(R_2\le r\le R_2+1)}\\
 & \les &  \| Z^{\le k} \pa^{\le 1} u(0, \cdot)\|_{L^2}+\|Z^{\le k} F \|_{L^1_t L^2_x}+\|\pa^{\le k}(1-\psi_{R_2+1})  u
\|_{ LE_{g}}\ ,
\end{eqnarray*}
where we have used the assumption $R_2\ge R_1+1$ in the last inequality.

To complete the proof, it remains to control
$\|\pa^{\le k}(1-\psi_{R_2+1})  u\|_{ LE_{g}}$ and $\|Z^{\le k}(1-\psi_{R_2})  u
\|_{ LE_{g}}$, which are trivially controlled  by
$\|\pa^{\le k} u\|_{ LE_{g}}$ and so is \eqref{eq-LE-2} by applying \eqref{H2}.
\end{prf}

As byproduct, let us record the fact that \eqref{H2} could be
replaced by even weaker assumptions, for example \eqref{H2'}.
\begin{prop}
\label{thm-H2'}
Let $n\ge 3$, then \eqref{H1} and \eqref{H2'} imply \eqref{H2}.
\end{prop}
\begin{prf}
As in Theorem \ref{thm-LE},  by \eqref{H1} and  Lemma \ref{thm-LE0},
with large enough $R_2$,
$$
\|\pa^{\le k} \psi_{R_2} u\|_{LE_g}\les
 \| \pa^{\le k+1} u(0, \cdot)\|_{L^2}+\|\pa^{\le k} F \|_{L^1_t L^2_x}+
 \|\pa^{\le k+1} u
\|_{ L^2_{t,x}(R_2\le r\le R_2+1)}
\ .$$
Then $\|\pa^{\le k} u\|_{LE_g}\le
\|\pa^{\le k} \psi_{R_2} u\|_{LE_g}+\|\pa^{\le k} (1-\psi_{R_2}) u\|_{LE_g}$,
which, by using the equation,
is controlled by $$
 \| \nabla^{\le k}\pa^{\le 1} u(0, \cdot)\|_{L^2}+\|\pa^{\le k} F \|_{L^1_t L^2_x}+
 \|\pa^{\le k}\pa u
\|_{ L^2_{t,x}(R_1\le r\le R_2+1)}
+ \|\pa^{\le k} u
\|_{ L^2_{t,x}( r\le R_2+1)}
\ .$$
The last two terms are bounded by the left of \eqref{H2'} and this completes the proof.
\end{prf}

Now we are ready to prove
Lemma \ref{thm-9-keystep1}.
At first, we observe that 
$$\|\<r\>^{-(n-1)/4} Z^{\le k} u\|_{L^2 L^2 L^2}\les
\|Z^{\le k} u\|_{LE_g}+\|\<r\>^{-(n-1)/4}\psi_R Z^{\le k} u\|_{L^2 L^2 L^2}\ ,$$
then in view of Theorem \ref{thm-LE} and $(n-1)/4>1/2$, it suffices for us to give the estimates of
\beeq\label{eq-WS-n=4}
 (\ln T)^{-1/q}\|\<r\>^{-1/q}\psi_R Z^{\le k} u\|_{L_T^{q}L^{2}}
+\|\<r\>^{-1/q-}\psi_R Z^{\le k} u\|_{L_T^{q}L^{2}}, q\ge 2  \ .\eneq

Let $w=\psi_R Z^{\le k} u$ with $R\ge R_3\ge R_2$ and $R_3$ to be fixed later, where 
$R_2$ is the constant appeared in the proof of Theorem \ref{thm-LE}.
Then, with $g_0=m (1-\psi_{R_3-1})+g \psi_{R_3-1}$,
$$\Box_{g_0} w=\Box_g w=[\Box_g,\psi_R] Z^{\le k}  u
+\psi_R [\Box_g, Z^{\le k}] u+\psi_R Z^{\le k} F\ .$$

For the first term on the right, we see that $[\Box_g,\psi_R] Z^{\le k}  u=\chi_{r\in [R,R+1]} \CO(\pa^{\le {k+1}} u)$.
Turning to the second term on the right, we note that
$$\psi_R [\Box_g, Z]u=
 r_1 \psi_R \pa^2  u+  r_2 \psi_R \pa  u\ ,
$$
where
\beeq\label{eq-lo}
\|r_1\|_{\ell^{1}_1 L^\infty_t L^\infty_x}+\||\pa r_1|+|r_2|\|_{\ell^{2}_1 L^\infty_t L^\infty_x}\les 1 .\eneq
Here, we used \eqref{H1.1} and the fact that the first perturbation of the metric is radial.
Hence, by induction,
$$\psi_R [\Box_g, Z^{\le j+1}]u=
 r_1 \psi_R \pa^2 Z^{\le j} u+  r_2 \psi_R \pa Z^{\le j} u\ ,\ \forall j\ge 0.
$$
Using the equation \eqref{lineqn}, when $k\ge 1$, we can get further
\beeq\label{eq-comm}
\psi_R [\Box_g, Z^{\le k}]u=
\nabla ( r_1 \psi_R \pa Z^{\le {k-1}} u)+  r_2 \psi_R \pa Z^{\le {k-1}} u+
r_1 \psi_R Z^{\le {k-1}} F
\ .
\eneq
In conclusion, we can write $\Box_{g_0} w$ as
$$\chi_{r\in [R,R+1]} \CO(\pa^{\le {k+1}} u)+
\nabla ( r_1 \psi_R \pa Z^{\le {k-1}} u)+  r_2 \psi_R \pa Z^{\le {k-1}} u+
r_1 \psi_R Z^{\le {k-1}} F+\psi_R Z^{\le k} F\ .$$

Then, by choosing $R_3$ large enough such that
$\Box_{g_0}$ satisfies the assumption in  Theorem \ref{thm-wStri-homo},
\eqref{eq-WS-n=4} could be bounded by
\begin{eqnarray*}
&&\|\psi_R Z^{\le k} u(0,\cdot)\|_{L^2}+\|\psi_R \pt Z^{\le k} u(0,\cdot)\|_{\dot H^{-1}}+\|\pa^{\le {k+1}} u\|_{L^2_T L^2_x(r\in [R, R+1])}
\\
 &  & 
+\|r_1 \psi_R \pa Z^{\le {k-1}} u\|_{\dot \ell^{1/2}_1 L^2_T L^{2}}
+\| r_2 \psi_R \pa Z^{\le {k-1}} u\|_{\dot \ell^{3/2}_2 L^2_T L^2_x}
\\
 &  &+\| (\psi_R-\psi_R^2)(r_1 Z^{\le {k-1}} F+Z^{\le k} F)\|_{ L^1_T \dot H^{-1}}\\
 &  &+\|r^{-(n-3)/2}(r_1 \psi_R^2 Z^{\le {k-1}} F+\psi_R^2 Z^{\le k} F)\|_{ L_T^2 L^1_r H^{-1/2+}_\omega}\\
&\les &\|\psi_R Z^{\le k} u(0,\cdot)\|_{L^2}+\|\psi_R \pt Z^{\le k} u(0,\cdot)\|_{\dot H^{-1}}+ \|\pa^{\le k} u\|_{LE_g}\\
  &  &
+
(\|r_1\|_{\ell^{1}_1 L^\infty_t L^\infty_x}
+\| r_2\|_{\ell^{2}_1 L^\infty_t L^\infty_x})
\|\psi_{R-1} \pa Z^{\le {k-1}} u\|_{ \ell^{-1/2}_\infty L^2_T L^{2}_x}\\
&&+
\| r Z^{\le k} F\|_{ L^1_T L^{2}(R\le r\le R+1)}
+\|r^{-(n-3)/2} \psi_R^2 Z^{\le k} F\|_{ L_T^2 L^1_r H^{-1/2+}_\omega}
\\
  & \les & 
  \|\psi_R Z^{\le k} u(0,\cdot)\|_{L^2}+\|\psi_R \pt Z^{\le k} u(0,\cdot)\|_{\dot H^{-1}}+  
  \|Z^{\le k} u\|_{LE_g}
+
\|F\|_{N^{k}_{R,p}}\ ,
\end{eqnarray*}
where we used the fact that $\psi_R-\psi_R^2$ is supported in $R\le r\le R+1$ and the duality to the Hardy inequality, $\|f\|_{\dot H^{-1}}\les \|rf\|_{L^2}$. 
The estimate for the third term on the right, $\|Z^{\le k} u\|_{LE_g}$,  is provided by Theorem \ref{thm-LE}.

Finally, we deal with the first two terms on the right in terms of the initial data.  By using the equation  \eqref{lineqn}, we have 
$\pt^{j+2}u(0)=\O(|\pa \nabla^{\le j+1}u(0)|+|\pa^{\le j} F(0)|)$, and
the first term is controlled by
\begin{eqnarray*} 
&&
 \|\psi_R Y^{\le k} u_0\|_{L^2}
 +\|\psi_R Y^{\le k-1} u_1\|_{L^2}
+\sum_{0\le j\le k-2} \|\psi_R Y^{\le k-j-2}\pt^{j+2} u(0,\cdot)\|_{L^2}\\
&
\les&
 \|\psi_R Y^{\le k} u_0\|_{L^2}
 +\|\psi_R Y^{\le k-1} u_1\|_{L^2}
+\sum_{0\le j\le k-2} \|\psi_R Y^{\le k-j-2}\pa^{\le j} F(0,\cdot)\|_{L^2}
\\
&
\les&
 \|\psi_R Y^{\le k} u_0\|_{L^2}
 +\|\psi_R Y^{\le k-1} u_1\|_{L^2}
+ \|\psi_R Z^{\le k-1} F(t,\cdot)\|_{L^1_T L^2}\ ,
\end{eqnarray*} where in the last step we used the Sobolev embedding $W^{1,1}\subset C^0\cap L^\infty$ for time.
Similarly, for the second term, it is bounded by
\begin{eqnarray*} 
&&
 \|\psi_R Y^{\le k} u_1\|_{\dot H^{-1}}
+\sum_{0\le j\le k-1} \|\psi_R Y^{\le k-j-1}\pt^{j+2} u(0,\cdot)\|_{\dot H^{-1}}\\
&\les
&
 \| Y^{\le k} u_0\|_{L^2}
+ \|\psi_R Y^{\le k} u_1\|_{\dot H^{-1}}
+\sum_{0\le j\le k-1} \|\psi_R Y^{\le k-j-1}\pt^{\le j} F(0,\cdot)\|_{\dot H^{-1}}\\
&\les
&
 \| Y^{\le k} u_0\|_{L^2}
+ \|\psi_R Y^{\le k} u_1\|_{\dot H^{-1}}
+ \|\psi_R Z^{\le k-1} F(0,\cdot)\|_{\dot H^{-1}}\ .
\end{eqnarray*}
This completes the proof.

\section{The space-times with stationary split metric near infinity}\label{sec-SplitMetric}
In this section, we give the proof of Theorem \ref{metaTheorem2}, for the space-times with stationary and split metric near spatial infinity. 

With more notations at hand, we could state a more precise version as follows
\begin{thm}\label{metaTheorem2.2} Let $n=3$, $2\le p\le 1+\sqrt{2}=p_c$, $s_d=1/2-1/p$,  and assume \eqref{H2} and \eqref{H1} with \eqref{eq-StaSplit}. For the problem \eqref{eq-nlw}, there exists $R, c, \ep_0>0$ such that 
 for any $(u_0, u_1)$ with 
 \beeq\label{eq-Sc-data}
 E[u_0, u_1] := \|Y^{\le 2} (u_0, u_1)\|_{H^1\times L^2} 
 + \|\psi_R Y^{\le 2}  u_1\|_{\dot H^{s_d-1}}=\ep\le \ep_0\ ,
\eneq
 there is a unique solution $u$ in $M\cap ([0, T_\ep]\times \R^3)$, 
 such that
$\| \psi_R Z^{\le 2} u \|_{WS_{p,\al}} +  \| Z^{\le 2} u\|_{LE_g}\le C_\al \ep$, where 
$T_\ep$ is defined in \eqref{eq-life1}, and for $T\ge 2$, 
\beeq\label{eq-def-WSpalpha}\|u\|_{WS_{p,\al}} \equiv 
\left\{
\begin{array}{ ll}
T^{-\al} \|r^{\frac{n-1}2-\frac np+\al} \psi_{1} u\|_{L^p_T L^p_r L^2_\omega}, &\al\in (0, 1/p],\\
(\ln T)^{-1/p} \|r^{\frac{n-1}2-\frac np} \psi_{1} u\|_{L^p_T L^p_r L^2_\omega}, &\al=0\ ,\\
\|r^{\frac{n-1}2-\frac np+\al} \psi_{1} u\|_{L^p_T L^p_r L^2_\omega}, &\al<0\ .
\end{array}
\right.\eneq
 \end{thm}
 
As in Section \ref{sec-General}, 
the novel part in the proof is the following version of the space-time estimates for the D'Alembertian with small, split, metric perturbation.
\begin{thm}[Space-time estimates]\label{ap-thm-le-0wS}
Let $n\ge 3$, $2< p<\infty$, $\al\le 1/p$ and assume \eqref{H1} with \eqref{eq-StaSplit}, then there exists $K\ge \tilde R_0$ such that
\begin{eqnarray}
\|u\|_{WS_{p,\al}\cap L^\infty \dot H^{s_d}} 
& \les_\al & \|u_0\|_{\dot H^{s_d}}+ \|u_1\|_{\dot H^{s_d-1}}+\|r^{-(n-3)/2-s_d}F_{1}\|_{L^2L^1H^{-1/p+}_{\omega}}
\nonumber
\\
 &&+\|F_{2}\|_{\dot \ell^{1/2}_1 L^2(L^2\cap \dot H^{-1})+\dot \ell^{1+1/p}_2 L^2L^2}
\ ,
\label{eq-WSp-Split}
\end{eqnarray}
for any $u$ satisfies $\Box_g u=F$ with initial data $(u_0, u_1)$ on $\{t=0\}$ and vanishes in the region $\{r<K\}$.
When $p=2$, we have instead
\beeq\label{eq-WSp-Split-9}
\|u\|_{WS_{2,\al}\cap L^\infty  L^2_x} 
\les_\al  \|w_0\|_{L^2}+ \|w_1\|_{\dot H^{-1}} 
+\|r^{-\frac{n-2}{2}} F_1\|_{L^1 L^1 H^{-\frac{1}{2}}_\omega}
+\|F_2\|_{\dot \ell^{1/2}_1 L^2 \dot H^{-1}+\dot \ell^{3/2}_1 L^2 L^2}\ .\eneq
\end{thm}
We note that the left hand side resembles the weighted Strichartz estimates \cite{FaWa, HMSSZ} (see also \cite{Hi07} for previous radial version).

\subsection{Proof of Theorem \ref{ap-thm-le-0wS}}
In this section, we give the proof of Theorem \ref{ap-thm-le-0wS}, which is divided into three parts. We first reduce the proof for operators of the form $\pt^{2}+P^{2}$, then deal with the homogeneous estimates and the inhomogeneous estimates separately.
\subsubsection{Theorem \ref{ap-thm-le-0wS}: reduction}
As in \cite{SoWa10} and \cite{LMSTW}, we want to rewrite the equation 
$\Box_g u=0$ near infinity as $(\pt^2+P^2) w=0$ so that $P$ is non-negative, elliptic and self-adjoint operator with respect to $L^2(\R^n)$.

Let $\tilde g=(1-\psi_{K-1} )m+\psi_{K-1} g$ with $K\ge \tilde R_0+1$ to be determined.
Since we are assuming that $u$ vanishes in the region $\{r<K\}$,
we have $\Box_g u=\Box_{\tilde g} u$.

Recall the condition \eqref{eq-StaSplit}
and \eqref{H1.1}, $\tilde g$ is stationary and split metric. 
Moreover, there exists $K_1\ge \tilde R_0+1$, such that
for  any $r\ge K_1$, we have $g_0(x)=-\tilde g_{00}\in (1/2, 2)$, $h_0(x)=
\det(\tilde g_{ij})\in (1/2, 2)$.
When $K\ge K_1$,
a direct calculation yields
$$
(\pt^2+P_K^2) w=- h_0^{1/4}g_0^{3/4} \Box_{\tilde g} u=- h_0^{1/4}g_0^
{3/4} \Box_{g} u$$
where $w=(h_0/g_0)^{1/4} u$ and $P_K$ is self-adjoint (with respect to $L^2(\R^n)$) and non-negative  with
$$P^2_K=-(g_0/h_0)^{1/4}\pa_i (g_0 h_0)^{1/2} \tilde g^{ij} \pa_j (g_0/h_0)^{1/4}\ .$$

Then, as in \cite{LMSTW}, the proof of Theorem \ref{ap-thm-le-0wS} is reduced to the following
 proposition.
\begin{prop}\label{ap-thm-le-0wS-1}
Let $n\ge 3$,
$T\ge 2$, $p\in [2,\infty)$ and $\al\le 1/p$.
Then there is a large enough $K>0$, such that we have
\eqref{eq-WSp-Split} and \eqref{eq-WSp-Split-9},
for any solutions $u$ to the equation $(\pt^2 +P_K^2)u=F$ with initial data $(u_0, u_1)$ on $\{t=0\}$.
\end{prop}

\subsubsection{Proof of Proposition \ref{ap-thm-le-0wS-1}: homogeneous estimates}

At first, we present a lemma, which asserts that $P_K$ behaves like $\nabla$ in appropriate function spaces. The proof is similar to that of \cite{LMSTW} and we omit it.
\begin{lem}\label{ap-thm-relate} For $P_K$, we have 
\beeq\label{ap-eq-equiv}\|P_K f\|_{\dot H^{s-1}}\les \|f\|_{\dot H^s},
\|P_K^{-1} f\|_{\dot H^{s}}\les \|f\|_{\dot H^{s-1}},
 s\in [0 , 1]\ ,\eneq
and
\beeq\label{ap-eq-0}\|P_K f\|_{\ell^{-1/2}_\infty L^2_{t,x}}
\les\|\nabla_x f\|_{\ell^{-1/2}_\infty L^2_{t,x}}+
\|f\|_{\ell^{-3/2}_\infty L^2_{t,x}}\ .\eneq
\end{lem}

Then, we notice that, provided that $K$ is large enough,  $\pt^2+P_K^2$ is of the form \eqref{eq-Box-per} with
\eqref{eq-kss-as-inh}, and so
Lemma \ref{thm-LE0}  gives us
$$
(\pt^2+P_K^2)u=0\Rightarrow
\|\pa u\|_{L^\infty_t L^2_x\cap \dot \ell^{-1/2}_\infty L^2_{t,x}}
+\| u\|_{\dot \ell^{-3/2}_\infty L^2_{t,x}}
\les \|\pa u(0,\cdot)\|_{L^2}\ .$$
Recall that $u=\cos (t P_K) u_0+P_K^{-1}\sin (t P_K) u_1$, with the help of Lemma \ref{ap-thm-relate}, we could extract the following estimates
\beeq\label{ap-eq-1}\|e^{it P_K} f\|_{L^\infty_t L^2_x\cap 
\ell^{-1/2}_\infty L^2_{t,x}}\les \|f\|_{L^2}, \eneq
This gives us the homogeneous part of \eqref{eq-WSp-Split-9} in light of \eqref{KSS} and \eqref{ap-eq-equiv}. In particular, for $T\ge 2$ and $\al\in (0,1/2]$,
\beeq\label{ap-eq-WS2}T^{-\al}\|\<r\>^{-\frac 12+\al} e^{it P_K} f\|_{L^2_T L^2_x}
+
(\ln T)^{-1/2}\|\<r\>^{-\frac 12} e^{it P_K} f\|_{L^2_T L^2_x}
+\|\<r\>^{-\frac 12-} e^{it P_K} f\|_{L^2_T L^2_x}
\les \|f\|_{L^2}. \eneq

Applying \eqref{ap-eq-equiv}, we get
$$\|e^{it P_K} f\|_{L^\infty_t \dot H^1_x}
=\|P_K^{-1}P_K e^{it P_K} f\|_{L^\infty_t \dot H^1_x}
\les \| e^{it P_K}P_K f\|_{L^\infty_t  L^2_x}
\les\|  P_K f\|_{ L^2_x}
\les \|f\|_{\dot H^1}, $$
and so 
\beeq\label{ap-eq-2}\|r^{(n-1)/2} e^{it P_K} f\|_{L^\infty_t L^\infty_r L^2_\omega}\les
\|e^{it P_K} f\|_{L^\infty_t \dot B^{1/2}_{2,1}}\les \|f\|_{\dot  B^{1/2}_{2,1}}, \eneq
by real interpolation and Lemma \ref{thm-trace}.

By interpolation between \eqref{ap-eq-WS2} and \eqref{ap-eq-2},  we get
\beeq\label{ap-eq-3}\|e^{it P_K} f\|_{WS_{p,\al}}\les \|f\|_{\dot  B^{1/2-1/p}_{2,1}}, 2<p<\infty, \al\le 1/p\ .\eneq
An application of real interpolation further (with $(\theta, p)$) gives us
\beeq\label{ap-eq-4}\|e^{it P_K} f\|_{WS_{p,\al}}\les \|f\|_{\dot  B^{1/2-1/p}_{2,p}}
\les \|f\|_{\dot  H^{1/2-1/p}}
, 2<p<\infty, \al\le  1/p\ .\eneq
This completes the proof of the homogeneous part of Proposition \ref{ap-thm-le-0wS-1},
by \eqref{ap-eq-equiv}.

\subsubsection{Proof of Proposition \ref{ap-thm-le-0wS-1}: inhomogeneous estimates}

At first, we apply the local energy estimates \eqref{eq-kss-inh}, backward in time,  for inhomogeneous wave equation $(\pt^2+P_K^2) w=G$ with vanishing data at $t=T$, we get
$$
\|w\|_{L^\infty \dot H^1\cap \dot\ell^{-1/2}_\infty L^2 \dot H^1
\cap \dot \ell^{-3/2}_\infty L^2 L^2}\les \|G\|_{L^1L^2+\dot \ell^{1/2}_1 L^2 L^2}\ ,$$
which, by duality, yields the following estimates 
\beeq\label{eq-Ws-05}
\|u\|_{WS_{2,\al}}\les \|u\|_{L^\infty L^2 \cap \dot \ell^{-1/2}_\infty L^2 L^2}\les
 \|F\|_{L^1\dot H^{-1}+\dot \ell^{1/2}_1 L^2 \dot H^{-1}
+\dot \ell^{3/2}_1 L^2 L^2}\ ,\ \al\le1/2
\eneq
for  $(\pt^2+P_K^2) w=G$  with vanishing data at $t=0$, which, combined with 
the duality to \eqref{eq-trace0}, gives us the inhomogeneous part of \eqref{eq-WSp-Split-9}.

For the case $p\in (2,\infty)$,
we claim that 
\beeq\label{ap-eq-6}\|e^{it P_K} f\|_{L^\infty \dot H^s\cap X^s}\les \|f\|_{\dot  H^{s}}, s\in (-1, 1], n\ge 3\ .\eneq
Based on \eqref{ap-eq-4}, \eqref{ap-eq-equiv}, duality of \eqref{ap-eq-6}, we have
\begin{eqnarray*}
\|\int_\R P_K^{-1}e^{i(t-\tau) P_K} F(\tau)d\tau\|_{WS_{p,\al}}
& \les & 
\|\int_\R P_K^{-1}e^{-i\tau P_K} F(\tau)d\tau\|_{\dot H^{s_d}}
 \\
& \les & \|\int_\R e^{-i\tau P_K} F(\tau)d\tau\|_{\dot H^{s_d-1}}\\
& \les & \|F\|_{(X^{1-s_d})'}\ .
\end{eqnarray*}
Recall that
Lemmas \ref{thm-MetTa29} and
\ref{thm-LE-key-observ}
 give us
$$\|u\|_{\dot \ell^{-1/2}_\infty L^2 (L^2+\dot H^{1})}
+\|u\|_{\dot \ell^{-1-1/p}_2 L^2  L^2}
+\|r^{(n-3)/2+s_d}u\|_{L^2L^\infty H^{1/p-}_{\omega}}
\les \|u\|_{X^{1-s_d}}\ .$$
By duality, we get 
$$\|\int_\R P_K^{-1}e^{i(t-\tau) P_K} F(\tau)d\tau\|_{WS_{p,\al}}
\les 
\|r^{-\frac{n-3}{2}-s_d}F_{1}\|_{L^2L^1H^{-\frac{1}{p}+}_{\omega}}
+\|F_{2}\|_{\dot \ell^{\frac{1}{2}}_1 L^2(L^2\cap \dot H^{-1})+\dot \ell^{1+\frac{1}{p}}_2 L^2L^2}\ .
$$
By Christ-Kiselev's lemma 
\cite{CK01},
this gives us the inhomogeneous part
of \eqref{eq-WSp-Split}.

It remains to prove the claim \eqref{ap-eq-6}.
Recall that for $(\pt^{2}+P_{K}^{2})u=0$, Lemma \ref{thm-MetTa9} gives us
$$\|\pa u\|_{L^\infty_t \dot H^s \cap X^{s}}\les \|u_0\|_{\dot H^{s+1}}+\|u_1\|_{\dot H^s}\ ,\ |s|<1 .$$
Recall that $u=\cos (t P_K) u_0+P_K^{-1}\sin (t P_K) u_1$, with the help of Lemmas \ref{ap-thm-relate} and \ref{thm-MetTa29}, we could extract the following estimates
$$\|\cos(t P_K) f\|_{L^\infty_t \dot H^{s} \cap X^{s}}
\les \|\nabla \cos(t P_K) f\|_{L^\infty_t \dot H^{s-1} \cap X^{s-1}}
 \les \|f\|_{\dot H^{s}}, s\in (0,2), $$
$$\|\sin(t P_K) f\|_{L^\infty_t \dot H^s \cap X^{s}}
\les \|\nabla P_{K}^{-1} \sin(t P_K) P_{K} f\|_{L^\infty_t \dot H^{s-1} \cap X^{s-1}}
 \les \|P_{K} f\|_{\dot H^{s-1}}
  \les \|f\|_{\dot H^s}, s\in (0, 1],$$
$$\|\cos(t P_K) f\|_{L^\infty_t \dot H^s \cap X^{s}}
=\|\pt P_{K}^{-1}\sin(t P_K) f\|_{L^\infty_t \dot H^s \cap X^{s}} 
 \les \|f\|_{\dot H^s}, |s|<1,$$
$$\|\sin(t P_K) f\|_{L^\infty_t \dot H^s \cap X^{s}}
=\|\pt \cos(t P_K) P_{K}^{-1} f\|_{L^\infty_t \dot H^s \cap X^{s}} 
 \les \|P_{K}^{-1} f\|_{\dot H^{s+1}}
  \les \|f\|_{\dot H^s}, s\in (-1, 0],$$
and so is \eqref{ap-eq-6}.

\subsection{Sketch of the proof for Theorem \ref{metaTheorem2.2}}
The proof proceeds as in Section \ref{sec-General}, with minor modifications. Here we give the sketch of the proof for $p\in (2, 1+\sqrt{2}]$ at first.
 Let $T>2$ and $\al\in [0, 1/p]$ to be determined, 
\begin{equation}
   \|\phi\|_{X_R^m} =  
  \| \psi_R Z^{\le m} \phi \|_{WS_{p,\al}\cap L^\infty \dot H^{s_d}}
  +
  \| Z^{\le m} \phi\|_{LE_g} 
+  
\|r^{-\frac 1p} \psi_R Z^{\le m} \phi\|_{L^p_T L^p_r L^2_\omega}
  ,\label{eq-10-X-norm}
\end{equation}
\beeq\label{eq-10-N-norm}
 \|F\|_{N^m} =  
 \| \psi_R^p r^{-s_d}Z^{\le m} F\|_{L^2L^1H^{-1/p+}_{\omega}}+\|Z^{\le m} F\|_{L^1L^2 L^2} \ .
\eneq
Notice that by \eqref{eq-trace2-sd} and H\"older's inequality, we have
\begin{eqnarray}
&&\|r^{1-5/(2p)+\al/2} \psi_R Z^{\le m} \phi \|_{L^{2p}_T L^{p} L^{2}}  \nonumber\\
&   \les&
(\|r^{1-3/p+\al} \psi_R Z^{\le m} \phi \|_{L^{p}_T L^{p} L^{2}} )^{1/2}
(\|r^{1-2/p} \psi_R Z^{\le m} \phi \|_{L^{\infty}_T L^{p} L^{2}} )^{1/2}\nonumber\\
&\les&
B_{T}\|\psi_R Z^{\le m} \phi\|_{WS_{p,\al}\cap L^\infty \dot H^{s_d}}\les   B_{T}    \|\phi\|_{X_R^m}, \al\in [0, 1/p],\label{eq-10-rel}
\end{eqnarray} 
here,
$$B_T=\left\{
\begin{array}{ll}
 T^{\al/2}     &   \al\in (0,1/p]\\
    (\ln T)^{1/(2p)}  &   \al=0\ .
\end{array}
\right.$$
Similar to the proof of Lemma \ref{thm-9-keystep1}, we get from Theorem \ref{thm-LE} and Theorem \ref{ap-thm-le-0wS} that there exists $R>R_1$ such that 
\begin{equation}\label{eq-10-linest}
\|u\|_{X_R^m} \lesssim \|Y^{\le m} u_0\|_{H^{1}}+\|Y^{\le m}u_1\|_{L^2}
+\|\psi_R Y^{\le m}u_1\|_{\dot H^{s_d-1}}
+\|\psi_R Z^{\le m-1}F(0,\cdot)\|_{\dot H^{s_{d}-1}}
+ \|F\|_{N^m} \ ,
\end{equation} for any $m\ge 0$ and solutions to
\eqref{lineqn}.
Moreover, we have
the nonlinear estimates
\begin{lem}\label{thm-10-keystep}
Let  $n=3$,  $p\in (2, 1+\sqrt{2}]$ and
$\al=2\frac{1+2 p-p^2}{p^2}$. Then
$\al\in [0, 1/p)$,
$p[1-5/(2p)+\al/2]=- s_d$,
\beeq
\|F_p(u)\|_{N^2}  \les   B_{T}^{p}  \|u\|_{{X}_{R}^2}^p \label{eq-10-keystep}\ ,\eneq
\beeq\label{eq-10-keystep2}
\|F_p(u)-F_p(v)\|_{N^{0}}\les 
B_{T}^{p} 
(\|u\|_{X_R^2}+\|v\|_{X_R^2})^{p-1}\|u-v\|_{X_R^0}\ ,\eneq
for any $T\ge 2$.
\end{lem}
Based on \eqref{eq-10-linest} and
Lemma \ref{thm-10-keystep}, through iteration in $X^k_R$, we can show existence of the solution with
$\|u\|_{X_R^2}\les \ep$, 
in the interval $t\in [0,
T_{\ep}]$,
provided that $B_{T_{\ep}}^{p} \ep^{p-1}\ll 1$. This gives us the long time existence up to
$T_{\ep}$ given in \eqref{eq-life1} and the proof of Theorem \ref{metaTheorem2.2} is reduced to that of  Lemma \ref{thm-10-keystep}.

For the remaining case $p=2$, we use \eqref{eq-WSp-Split-9} instead of \eqref{eq-WSp-Split}. More precisely, let $T>2$ to be determined, 
 we introduce
\begin{equation}
   \|\phi\|_{X_R^m} =  
     \| \psi_R Z^{\le m} \phi \|_{WS_{2,1/4}\cap
     WS_{2,0}
     \cap L^\infty L^2}
  +
  \| Z^{\le m} \phi\|_{LE_g} 
  ,
\end{equation}
\beeq
 \|F\|_{N^m} =  
 \| \psi_R^p r^{-1/2}Z^{\le m} F\|_{L^1L^1H^{-1/2}_{\omega}}+\|Z^{\le m} F\|_{L^1L^2 L^2} \ .
\eneq
Moreover, we have the following 
variant of Lemma \ref{thm-10-keystep},
\beeq
\|F_2(u)\|_{N^2} 
\les T^{1/2} \|u\|_{{X}_{R}^2}^2 \ ,\eneq
\beeq\label{eq-10-keystep2}
\|F_2(u)-F_2(v)\|_{N^{0}}\les 
T^{1/2} 
(\|u\|_{X_R^2}+\|v\|_{X_R^2}) \|u-v\|_{X_R^0}\ ,\eneq
for any $T\ge 2$. Equipped with these estimates, 
we can show existence of the solution with
$\|u\|_{X_R^2}\les \ep$, 
in the interval $t\in [0,
T_{\ep}]$,
provided that $T_{\ep}^{1/2} \ep^{2}\ll \ep$, which completes the proof of Theorem \ref{metaTheorem2.2}.

\subsection{Proof of Lemma \ref{thm-10-keystep}}
In this subsection, we give the proof of \eqref{eq-10-keystep}, and omit the similar proof of \eqref{eq-10-keystep2}.

Let $\be=5/(2p)-\al/2-1$, 
 we  observe that
 $\be p=s_d$ and
$$|Z^{\le 2} F_p(u)|\lesssim  |u|^{p-1} |Z^{\le 2} u|
+
|u|^{p-2} |Z^{\le 1} u|^{2}\equiv I+II
\ .
$$ 

We start with bounding the first term in \eqref{eq-10-N-norm}.  We first
note that
\begin{eqnarray*}
 \| r^{-\be p} \psi_R^p Z^{\le 2}
F_p(u)\|_{L^2L^1 L^2} 
 & \les & \|r^{-\be} \psi_R u\|^{p-1}_{L^{2p}L^pL^\infty}
\|r^{-\be} \psi_R Z^{\le 2} u\|_{L^{2p}L^pL^2}
\\
&&+\|r^{-\be} \psi_R u\|^{p-2}_{L^{2p} L^pL^\infty}
\|r^{-\be} \psi_R Z^{\le 1} u\|_{L^{2p} L^pL^4}^{2}.
\end{eqnarray*}
By \eqref{eq-10-rel}, the Sobolev embeddings on the unit sphere $\Sp^{2}$, $H^2_\omega\subset L^\infty_\omega$ and $H^1_\omega\subset
L^4_\omega$, it follows that 
\beeq\label{eq-10-bd}
 \| \psi_R^p r^{-s_d}Z^{\le m} F_p(u)\|_{L^2L^1H^{-1/p+}_{\omega}}\les
B_T^p \|u\|_{X_{R}^{2}}^{p}\ .\eneq

We now proceed to the second term in \eqref{eq-10-N-norm}.  We first observe
that
\begin{eqnarray*}\|I \|_{L^1L^2L^2(r\ge R+2)}&
  \lesssim &\|r^{\frac{1}{p(p-1)}}
u\|^{p-1}_{L^pL^{\frac{2p(p-1)}{p-2}}L^\infty(r\ge R+2)}
 \|r^{-1/p} \psi_R Z^{\le 2} u\|_{L^pL^pL^{2}}\\
 &\les&
 \|r^{-1/p} \psi_R Z^{\le 2} u\|_{L^pL^pL^{2}}^p
 \les \|u\|_{X_{R}^2}^{p}\ ,
\end{eqnarray*}
where we have applyed \eqref{eq-Sobo} in the second to the last inequality.
Moreover, by Sobolev embedding $\dot H^{2}\cap \dot H^1\subset L^\infty$,
\begin{eqnarray*}\|I\|_{L^1L^2L^2(r\le
  R+2)} &\lesssim&
\|u\|_{(L^2\cap L^{\infty})L^\infty L^\infty(r\le R+2)}^{p-1}
\|\partial^{\le 2} u\|_{L^2 L^2L^2(r\le R+2)}\\
&\les &\|u\|_{X_{R}^2}^{p}.
\end{eqnarray*}

Similarly, let $m_2=-\frac 1p+\frac{n-1}p-\frac{n-1}4$, we have
$(-\frac 1p+\frac{n-1}p)(p-2)+2 m_2=0$.
Thus by Lemma \ref{thm-Sobo}, $\|II\|_{L^1L^2L^2(r\ge
  R+2)}$ is bounded by
$$ \|r^{-\frac 1p+\frac{n-1}p} u\|^{p-2}_{L^p L^\infty
    L^\infty(r\ge R+2)}
\|r^{m_2} Z^{\le k-1} u\|_{L^p L^4 L^4(r\ge R+2)}^{2}
\les \|r^{-1/p} \psi_R Z^{\le 2} u\|_{L^pL^pL^{2}}^p
 \les \|u\|_{X_{R}^2}^{p}
.$$
It remains to bound
$$\|II\|_{L^1L^2L^2(r\le
  R+2)} \les \|u\|^{p-2}_{L^\infty L^\infty
    L^\infty(r\le R+2)}\|Z^{\le 1} u\|_{L^2
    L^4 L^4(r\le R+2)}^{2}\ .
$$ Using Sobolev embeddings $\dot H^1\cap \dot H^2\subset L^\infty$ and $H^1\subset L^4$, this is estimated by  $\|u\|_{X_{R}^2}^{p}$.
This completes the proof of Lemma \ref{thm-10-keystep}.

\section{Long time existence with small, compact supported data}
\label{sec-compact}

In this section, 
we prove Theorem \ref{metaTheorem2}, for the case when \eqref{eq-StaSplit} is not satisfied.  As before, we state a more precise version of
the theorem as follows.
\begin{thm}\label{metaTheorem}
 Let $n=3$, $2\le p\le 1+\sqrt{2}$, $\al=(1+2 p-p^2)/p^{2}\in [0,1/p]$, $R\ge R_3$
and assume \eqref{H1}, \eqref{H2}. Consider the problem \eqref{eq-nlw} with initial data
$(u_0, u_1)\in H^3\times H^2$ with compact support in $B_R$.
Then there exists $C>c>0$ and $\ep_0>0$,  such that
for any $(u_0, u_1)$ with \beeq\label{Th-as}\|u_0\|_{H^3}+\|u_1\|_{H^2}=\ep\le \ep_0\ ,\eneq
 there is a unique solution $u$ in $M\cap ([0, T_\ep]\times \R^3)$, such that
$\| \psi_R Z^{\le 2} u \|_{WS_{p,\al}} +  \| Z^{\le 2} u\|_{LE_g}\le C\ep$, here $T_\ep$ is defined in \eqref{eq-life2}.
\end{thm}

As in the previous sections,
the result relies on the following space-time estimates for the solutions of the inhomogeneous wave equation with small metric perturbation.
\begin{thm}
\label{thm-wStri}
Let $n\ge 3$,  $\Box_h$ as in \eqref{eq-Box-per1},
 with $h^{\al\be}_0=h_0^{\be\al}\in C^2$, $h_1\in C^{1}$. Then 
there exists a $\delta_1>0$, such that if $h$ satisfies
\eqref{eq-kss-as2'}, we have
\beeq\label{eq-WS2}\|u\|_{L^\infty_t L^2_x\cap \dot \ell^{-1/2}_\infty L^2_t L^2_x}\les
\|F\|_{L^1_t \dot H^{-1}+\dot \ell^{1/2}_1 L^2 \dot H^{-1}+\dot \ell^{3/2}_1 L^2 L^2}\ ,\eneq
\beeq
\|u\|_{WS_{p,\al}} \les\|F\|_{\dot \ell^{1/2}_1 L^2_T (\dot H^{-1}\cap L^{2})+L^1_T \dot B^{-1/2-1/p}_{2,1}
+\dot \ell^{1+1/p}_1 L^2_T L^2_x}\ ,
\label{eq-WSp}
\eneq
for any $T\ge 2$, $p\in (2,\infty)$, $\al\le 1/p$ and
solutions to $\Box_h u=F$ with vanishing initial data.
\end{thm}

Combined with 
Lemma \ref{thm-trace}, 
as a corollary of Theorem \ref{thm-wStri}, we have the following version of the inhomogeneous weighted Strichartz estimates
\begin{coro}
\label{thm-wStri2}
Under the same assumption in Theorem \ref{thm-wStri}, with $\Box_h u=F+\nabla G$, we have
\beeq\label{eq-WSp2}
\|u\|_{WS_{p,\al}} \les\|F\|_{\dot \ell^{1/p-(n-1)/2}_1 L_T^1 L^1_r L^{2}_\omega
+\dot \ell^{1+1/p}_1 L^2_T L^2_x}+\|G\|_{\dot \ell^{1/2}_1 L^2_T H^{1}}
\ ,
\eneq
for any $p\in [2,\infty)$ and $\al\in [0, 1/p]$.
\end{coro}

\subsection{Sketch of the proof for Theorem \ref{metaTheorem}}
The proof proceeds as in Section \ref{sec-SplitMetric}, with the 
 Theorem \ref{ap-thm-le-0wS} replaced by
Corollary \ref{thm-wStri2}. Here we give a sketch of the proof for $p\in [2, 1+\sqrt{2}]$.
 Let $T>2$ to be determined,
  we introduce the solution space $X^{k}_{R,p,\al}$ in the current setting,
 \begin{equation}
   \|\phi\|_{X^{k}_{R,p,\al}} =   
  \| \psi_R Z^{\le k} \phi \|_{WS_{p,\al}} +
  \| Z^{\le k} \phi\|_{LE_g} ,\label{X-norm}
\end{equation}
\beeq\label{N-norm}
 \|g\|_{N_{R,p}^{k}} =  \| \psi_R^p Z^{\le k} g \|_{\dot \ell^{1/p-1}_1 L_T^1 L^1 L^{2}_\omega} + \|Z^{\le k} g\|_{L^1 L^2 L^2} \ .
\eneq

Similar to the proof of Lemma \ref{thm-9-keystep1}, by Theorem\ref{thm-LE} and Corollary \ref{thm-wStri2}, we get
\begin{lem}\label{thm-key}
Let $n= 3$, $k\ge 0$, $p\in [2,\infty)$, $\alpha\le 1/p$ and assume \eqref{H1} and \eqref{H2}. Then there exists a constant $R_3\ge R_2$, where 
$R_2$ is the constant appeared in the proof of Theorem \ref{thm-LE},
 such that for any $R\ge R_3$ and $T\ge 2$, we have
\begin{equation}\label{linest}
\|u\|_{X^{k}_{R,p,\al}} \lesssim \|u_0\|_{H^{k+1}}+\|u_1\|_{H^k}+ \|F\|_{N^{k}_{R,p}} \ .
\end{equation}
Here $u$ can be any solutions to
\eqref{lineqn}, where $u_0$, $u_1$, $\pt^{\le {k-1}} F(0,r\omega)$ are supported in the region $\{r\le  R\}$.
\end{lem}
For $p\in [2, 1+\sqrt{2}]$,
let $\al=\frac{1+2 p-p^2}{ p^2}\in [0, 1/p]$, $\be=\frac 3p-1-\al$,
and $$C_T=\left\{
\begin{array}{ll}
 T^\al  ,   &   \al\in (0, 1/p] ,\\
    (\ln T)^{1/p} ,  &   \al=0\ .
\end{array}
\right.$$   Here the choice of $\al$ ensures 
$(1-\frac 3p+\al)p=-\be p=\frac 1p-1$.
With these notations, we see that
\begin{equation}
   \|\phi\|_{X_{R,p,\al}^{k}} \ge 
  C_T^{-1}\| r^{-\be} \psi_R Z^{\le k} \phi 
  \|_{L^p_T L^p_r L^2_\omega}
+  \| Z^{\le k} \phi\|_{LE_g},\label{X-norm2}
\end{equation}
\beeq\label{N-norm2}
 \|g\|_{N^{k}_{R,p}}= \| r^{-\be p} \psi_R^p Z^{\le k} g \|_{ L_T^1 L^1 L^{2}_\omega} + \|Z^{\le k} g\|_{L^1 L^2 L^2}\ .
\eneq 

For the proof of Theorem \ref{metaTheorem} by iteration in $X^2_{R,p,\al}$, recalling Lemma \ref{thm-key}, it suffices to prove the following lemma.
\begin{lem}\label{thm-keystep}
Let  $n=3$, $p\in
[2, 1+\sqrt{2}]$ and $\al=\frac{1+2 p-p^2}{ p^2}$. Then
\beeq\label{eq-keystep}\|F_p(u)\|_{N^2_{R,p}}\les C_T^p \|u\|_{{X}^2_{R,p,\al}}^p\ ,\eneq
\beeq\label{eq-keystep2}\|F_p(u)-F_p(v)\|_{N^{0}_{R,p}}\les C_T^p(\|u\|_{{X}^2_{R,p,\al}}+\|v\|_{{X}^2_{R,p,\al}})^{p-1}\|u-v\|_{{X}^{0}_{R,p,\al}}\ .\eneq
\end{lem}
The proof of Lemma \ref{thm-keystep} follows a similar way as that of Lemma \ref{thm-10-keystep} and we leave it for readers.

With all this lemmas available, it is easy to see that as long as
\beeq\label{eq-life-2}C_T^p \ep^{p-1}\ll 1\eneq
is satisfied, we could prove existence and uniqueness of solutions in $X_{R,p,\al}^2$ in the interval $t\in [0, T]$. The condition is satisfied for $T_\ep$ given in \eqref{eq-life2} and
this completes the proof of Theorem \ref{metaTheorem}.

\subsection{Inhomogeneous weighted Strichartz estimates}\label{Sec-sub-WS-inh}
As in \cite{LMSTW}, we will prove Theorem \ref{thm-wStri}, based on local energy estimates, duality and interpolation.

At first, 
as the proof of \eqref{eq-Ws-05} for $\pt^2+P_K^2$,
we get from the local energy estimates \eqref{eq-kss-inh}
with duality that
\beeq\label{eq-Ws-0}
\|u\|_{WS_{2,\al}}\les \|u\|_{L^\infty L^2 \cap \dot \ell^{-1/2}_\infty L^2 L^2}\les
 \|F\|_{L^1\dot H^{-1}+\dot \ell^{1/2}_1 L^2 \dot H^{-1}
+\dot \ell^{3/2}_1 L^2 L^2}\ ,\ \al\le1/2\ ,
\eneq
\beeq\label{eq-Ws-1}
\|u\|_{L^\infty \dot H^1} \les \|F\|_{L^1L^2+\dot \ell^{1/2}_1 L^2 L^2}\ ,
\eneq
for $\Box_h u=F$ with vanishing data at $t=0$, which gives us \eqref{eq-WS2}.
Now real interpolation between
\eqref{eq-Ws-0} and \eqref{eq-Ws-1}
with $(\theta,q)=(1/2,1)$, together with the endpoint trace estimates \eqref{eq-trace0}, give us
\beeq\label{eq-Ws-half}
\|r^{(n-1)/2}u\|_{L^\infty L^\infty L^2}\les
\|u\|_{L^\infty \dot B^{1/2}_{2,1}}\les \|F\|_{L^1\dot B^{-1/2}_{2,1}+
\dot \ell^{1/2}_1 L^2 (\dot H^{-1}\cap L^2)
+\dot \ell^{1}_1 L^2 L^2}\ .
\eneq

By H\"older's inequality, combining with
\eqref{eq-Ws-0}, for $\al\le 1/p$ with $2<p<\infty$, we arrive at
$$
\|u\|_{WS_{p,\al}}\les
 \|F_1\|_{L^1\dot H^{-1}}^{2/p} \|F_1\|_{L^1\dot B^{-1/2}_{2,1}}^{1-2/p}
  +
\|F_2\|_{\dot\ell^{1/2}_1 L^2 (\dot H^{-1}\cap L^2)}
+\|F_3\|_{\dot\ell^{3/2}_1 L^2 L^2}^{2/p} \|F_3\|_{\dot\ell^{1}_1 L^2 L^2}^{1-2/p}\ ,
$$
for $\Box_h u=F_1+F_2+F_3$ with vanishing data at $t=0$. After summation with spatial or frequency localization, we could improve the estimates to
\beeq\label{eq-Ws-p}
\|u\|_{WS_{p,\al}}\les
 \|F_1\|_{L^1\dot B^{-1/2-1/p}_{2,1}}
+
\|F_2\|_{\dot\ell^{1/2}_1 L^2 (\dot H^{-1}\cap L^2)}
+\|F_3\|_{\dot\ell^{1+1/p}_1 L^2 L^2}\ ,
\eneq
which is exactly \eqref{eq-WSp}.

\section{Discussion}\label{sec-Discuss}
As we have mentioned in the introduction, our main theorems (Theorems \ref{metaTheorem1}-\ref{metaTheorem2}, as well as their precise version
Theorems  \ref{metaTheorem1-2}, 
\ref{metaTheorem2.2} and \ref{metaTheorem})
are conditional results and apply for many space-times of interest. Remarkable examples including Schwarzschild/Kerr black hole backgrounds, as well as their high dimensional counterparts, nontrapping exterior domain, nontrapping asymptotically Euclidean manifolds, and asymptotically flat manifolds with small time-dependent metric perturbation.

\subsection{Schwarzschild space-time}
As we know, Schwarzschild space-time is asymptotically flat, stationary and split near spatial infinity. Moreover, the hypothesis \eqref{H2} has been verified in \cite{BS, BS2}, \cite{DaRo, DaRo09}, \cite{MMTT}.
Then, we can apply Theorem \ref{metaTheorem2} for \eqref{eq-nlw} with general small data to obtain the following bound for the lifespan $L_\ep$:
$$L_\ep\ge \left\{
\begin{array}{ll }
 c\ep^{-\frac{p(p-1)}{1+2p-p^2}},   &   2\le p<1+\sqrt{2}, \\
\exp(c\ep^{-2\sqrt{2}}),     &   p=1+\sqrt{2} \ .
\end{array}
\right.$$
Currently, there are no any known explicit upper bound for the lifespan, except
a weak result of nonexistence of global existence in \cite{CG06}.
In contrast,
the global existence for $p>1+\sqrt{2}$ has been verified for the Schwarzschild space-time in \cite{LMSTW}, with previous works 
\cite{DaRo05, BlSt06}.

 In addition, for the high dimensional (hyperspherical) Schwarzschild backgrounds, the local energy assumption \eqref{H2} is known from the works \cite{LM, Schlue}, and so is the almost global existence, $L_\ep\ge \exp(c\ep^{-2})$, for \eqref{eq-nlw} with $p=2$ and $n=4$,
and global existence for \eqref{eq-nlw} with $p=2$ and $n\ge 5$,
  by Theorem \ref{metaTheorem1}.

\subsection{Kerr space-times}
Turning to Kerr space-times, we know it is stationary, satisfying \eqref{H1}.
Concerning \eqref{H2}, it has been shown for slowly rotating Kerr space-times in \cite{TT}, \cite{AB09}, \cite{DaRo08, DaRoNotes, DaRoNew}. See \cite{TT} for the higher order ($k>0$) case.
Here, we also note that \cite{MTT} extend \eqref{H2} for a class of small perturbations of Kerr metrics, and  \cite{DaRoSR} essentially verifies \eqref{H2} for the full subextremal case $|a|<M$.

As a consequence, Theorem \ref{metaTheorem2} with lower bound
\eqref{eq-life2}
applies for the slowly rotating Kerr space-times,  concerning \eqref{eq-nlw} with $2\le p\le 1+\sqrt{2}$.
On the other hand, 
the global existence
with $p>1+\sqrt{2}$ is known in \cite{LMSTW, MW}.
It is also related to the open problem 11 in \cite[page 173]{DaRoNotes}, with $V(\psi)=\pm |\psi|^{p+1}$.

In addition, as generalization of Kerr backgrounds in higher dimension, \eqref{H2} for $(1+4)$-dimensional Myers-Perry space-times with small angular momenta is essentially verified in \cite{LMTT}. Thus, in this setting, we have almost global existence, $L_\ep\ge \exp(c\ep^{-2})$, for \eqref{eq-nlw} with $p=2$ and $n=4$, by Theorem \ref{metaTheorem1}.

\subsection{Nontrapping exterior domain} 
For the problem in the nontrapping exterior domain (even with weak trapping), with homogeneous Dirichlet boundary condition, the local energy assumption is well-known, see, e.g., \cite{HMSSZ, Yu11} for history in detail.
It is known that we have global existence \cite{DMSZ, HMSSZ, SmSoWa12} for $p>p_c$ and $2\le n\le 4$, when the exterior domain is nontrapping with compact metric perturbation.

Concerning the nonexistence of global solutions and the upper bound of the lifespan, for $g=m$, it is known that
$$L_\ep< \left\{
\begin{array}{ll }
 C \ep^{-\frac{2p(p-1)}{2+(n+1)p-(n-1)p^2}},   &   p<p_c(n), n\ge 3, 
\textrm{  \cite{ZhouHan1}}
 \\
\exp(C\ep^{-2-\sqrt{2}}),     &   p=1+\sqrt{2}, n=3, \textrm{\cite{LaZh14p}} \ .
\end{array}
\right.$$
On the other hand, for star-shaped obstacles, the following lower bounds are known 
$$L_\ep\ge  \left\{
\begin{array}{ll }
 c \ep^{-2},   &   p=2, n= 3, 
\textrm{  \cite{DuZh08}}
 \\
 c \ep^{-\frac{p(p-1)}{1+2p-p^2}},   &   2<p< 1+\sqrt{2}, n=3, 
\textrm{  \cite{Yu11}}
 \\
\exp(c\ep^{-\sqrt{2}}),     &    p=1+\sqrt{2}, n=3,  \textrm{\cite{Yu11}} 
 \\
\exp(c\ep^{-2}),     &   p=2, n=4, F_p=u^2, \textrm{\cite{ZhZh15}} 
 \\
\infty,     &   p=2, n\ge 5, F_p=u^2, \textrm{\cite{MetSo10}}
\end{array}
\right.$$
which agree with the upper bound for $2\le p<p_c(3)$.

Our results recover and strengthen the known existence results to any nontrapping exterior domain
\beeq\label{life}L_\ep\ge  \left\{
\begin{array}{ll }
 c \ep^{-\frac{p(p-1)}{1+2p-p^2}},   &   2\le p< 1+\sqrt{2}, n=3\ ,
 \\
\exp(c\ep^{-2\sqrt{2}}),     &    p=1+\sqrt{2}, n=3\ ,
 \\
\exp(c\ep^{-2}),     &   p=2, n=4, F_p=u^2\ ,
\\
\infty,     &   p=2, n\ge 5, F_p=u^2\ .
\end{array}
\right.\eneq

\subsection{Nontrapping asymptotically Euclidean manifolds}
Let $R_0=0$, $g=-dt^2+g_{ij}(x)dx^i dx^j$, and assume that there exists some $\rho>0$, such that
$$g_{ij}(x)=\delta_{ij}+g_{1,ij}(|x|)+g_{2,ij}(x), \nabla^a g_{k,ij}=\CO(\<r\>^{1-k-\rho-|a|})\ .$$
Moreover, assume the metric is nontrapping, we call such manifolds as nontrapping asymptotically Euclidean manifolds, which is of course stationary and split.

For such manifolds, it is known from \cite{BoHa} that the assumption \eqref{H2'}, with $R_1=0$ and $N>1/2$, is satisfied. By Proposition \ref{thm-H2'}, we have \eqref{H2}.

Then, we can apply Theorem \ref{metaTheorem1}-\ref{metaTheorem2} for \eqref{eq-nlw} with general small data to obtain the lower bound \eqref{life}.
Notice that it improves the corresponding results in \cite{WaYu11} in many ways.
First of all, the expected sharp lifespan is obtained for $2\le p<1+\sqrt{2}$ and $n=3$. Secondly, the metric assumption has been relaxed. Last, it is the first time to obtain the almost global existence for $p=p_c$ for $n=3,4$, 
in the setting of asymptotically Euclidean manifolds.

\subsection{Time-dependent, small, asymptotically flat manifolds}
In the case where $R_0=0$ and the constants in \eqref{H1.1} are sufficiently small,
the local energy hypothesis \eqref{H2} with $R_1=0$, follows from \eqref{eq-kss-inh}
and the same proof as Theorem \ref{thm-LE}.

As for the Kerr space-times, Theorem \ref{metaTheorem1}-\ref{metaTheorem2} apply in this setting, concerning \eqref{eq-nlw} with $2\le p\le p_c$.


\end{document}